\documentclass[12pt,a4paper]{amsart}

\usepackage{euscript,amsfonts,amssymb,amsmath,amscd,amsthm,enumerate,hyperref}

\usepackage{comment}

\usepackage{mathrsfs}

\usepackage{graphicx}

\usepackage{color}

\usepackage[margin=1.1in]{geometry}

\usepackage{bbm}

\newtheorem{theorem}{Theorem}[section]
\newtheorem*{theorem*}{Theorem}
\newtheorem{lemma}[theorem]{Lemma}
\newtheorem{definition}[theorem]{Definition}
\newtheorem{proposition}[theorem]{Proposition}

\newtheorem{assumption}[theorem]{Assumption}

\newtheorem{cor}{Corollary}
\newtheorem{ansatz}{Ansatz}
\newtheorem{ass}{Assumption}

\theoremstyle{remark}
\newtheorem{rmk}[theorem]{Remark}

\newcommand{\eps}{\varepsilon}

\newcommand{\E}{\mathbb E}
\newcommand{\EE}{\mathbb E}

\newcommand{\pp}{\mathbb{P}}

\newcommand{\rr}{\mathbb{R}}

\newcommand{\nn}{\mathbb{N}}
\newcommand{\ep}{\hfill \ensuremath{\Box}}
\newcommand{\eq}{\begin{equation}}
\newcommand{\en}{\end{equation}}

\newcommand{\bone}{\mathbf{1}}
\newcommand{\PP}{\mathbb{P}}
\newcommand{\RR}{\mathbb{R}}
\newcommand{\mX}{\mathcal{X}}
\newcommand{\mS}{\mathcal{S}}
\newcommand{\mD}{\mathcal{D}}
\newcommand{\bC}{\bar{C}}
\newcommand{\bnu}{\bar{\nu}}
\newcommand{\mV}{\mathcal{V}}
\newcommand{\hr}{\hat{r}}
\newcommand{\hb}{\hat{b}}
\newcommand{\blambda}{\bar{\lambda}}
\newcommand{\hC}{\hat{C}}
\newcommand{\hnu}{\hat{\nu}}

\numberwithin{equation}{section}

\title[Interaction through hitting times]{Mean Field Systems on Networks, with Singular Interaction Through Hitting Times} 

\author{Sergey Nadtochiy}
\address{Department of Applied Mathematics, Illinois Institute of Technology, Chicago, IL 60616, USA}
\email{snadtochiy@iit.edu}
\thanks{S. Nadtochiy is partially supported by the NSF CAREER Grant DMS-1855309.}

\author{Mykhaylo Shkolnikov}
\address{ORFE Department, Bendheim Center for Finance, and Program in Applied \& Computational Mathematics, Princeton University, Princeton, NJ 08544, USA}
\email{mshkolni@gmail.com}
\thanks{M. Shkolnikov is partially supported by the NSF grant DMS-1506290 and a Princeton SEAS innovation research grant.}


\begin{document}

\begin{abstract}
Building on the line of work \cite{Delarue1}, \cite{Delarue2}, \cite{NadShk}, \cite{DeTs}, \cite{HaLeSo}, \cite{HaSo} we continue the study of particle systems with singular interaction through hitting times. In contrast to the previous research, we (i) consider very general driving processes and interaction functions, (ii) allow for inhomogeneous connection structures, and (iii) analyze a game in which the particles determine their connections strategically. Hereby, we uncover two completely new phenomena.
First, we characterize the ``times of fragility" of such systems (e.g. the times when a macroscopic part of the population defaults or gets infected simultaneously, or when the neuron cells ``synchronize") explicitly in terms of the dynamics of the driving processes, the current distribution of the particles' values, and the topology of the underlying network (represented by its Perron-Frobenius eigenvalue).
Second, we use such systems to describe a dynamic credit-network game and show that, in equilibrium, the system regularizes, i.e. the times of fragility never occur, as the particles avoid them by adjusting their connections strategically.
Two auxiliary mathematical results, useful in their own right, are uncovered during our investigation: a generalization of Schauder's fixed-point theorem for the Skorokhod space with the M1 topology, and the application of max-plus algebra to the equilibrium version of the network flow problem.
\end{abstract}

\maketitle

\section{Introduction}
\label{se:intro}

In this paper, we continue our investigation of the mean field particle systems with singular interaction through hitting times.
The finite particle systems of such type are characterized by the following feature. Whenever a particle (represented by a stochastic process on the real line) hits a given barrier, it shifts the positions of all other particles instantaneously in the direction of the barrier. 
As a result, a single initial hitting time may cause an instantaneous ``cascade" of jumps. It turns out that, unlike in classical mean field systems, these cascades do not disappear as the population size increases to infinity, and they manifest themselves in the jumps of a representative particle. Such jumps are, often, connected to important real-world phenomena (e.g. phase transitions), hence, their presence in a model is necessary in order to analyze these phenomena. On the other hand, such discontinuities are the main source of the associated mathematical challenges.

The particle systems of such type arise in a variety of contexts. In \cite{Delarue1}, \cite{Delarue2}, they are used as neuron firing models with the purpose of capturing the synchronization phenomenon, namely, a macroscopic number of neurons in a part of the brain firing simultaneously. The latter, in turn, has been found to be correlated with many cognitive functions. Other particle systems interacting through hitting times are introduced in \cite{DeTs}, followed by \cite{DelNadShk}, and are motivated by the earlier research \cite{FaPr}, \cite{FaPr2}, \cite{FaPr3}, \cite{FaPr4} on supercooled liquids. Thereby, the goal is to describe the rapid advancement of the liquid-solid boundary when the supercooled liquid is located next to the warmer corresponding ice. In our article \cite{NadShk} and the papers \cite{HaLeSo}, \cite{HaSo} that followed, the interaction through hitting times represents the losses banks incur due to the defaults of other banks, which may cause further defaults and potentially trigger a default cascade of macroscopic size. The common feature of the three settings lies in the self-excitation on the microscopic (particle) level that, under appropriate circumstances, leads to discontinuities on the macroscopic level.

In addition to considering very general driving processes and interaction functions, the main contribution of the present work is two-fold. First and foremost, we allow for an inhomogeneous connection structure across particles. Recall that, in the mean field particle systems described in the previous paragraphs, the hitting time of any particle has a direct and equal effect on the positions of all other particles. This implies that any two particles are directly connected, and the strength of each connection is the same. Herein, on the contrary, we model the connections via an infinite directed weighted graph, where the weight of an edge indicates the strength of a connection. Specifying the weights of the graph, we can produce systems with various connection structures. This generalization allows one to produce more realistic models, e.g. it is clear that not all neurons are directly connected to each other, and that not all banks have equal exposures to each other. To the best of our knowledge, the particle systems of such form have not been considered in the existing literature.
It is worth mentioning that the connection structures used herein appear as the limits of finite random graphs in \cite{InhomRandGraphs}, and the particle systems with such connection structures are analyzed, e.g. in \cite{NagasawaTanaka1}, \cite{NeuroNonhom}, \cite{Wu2}, \cite{Wu1}, \cite{CSY}. However, the latter systems are of a classical, smooth, type, without the singular interaction through hitting times.
In Section \ref{se:existence}, we discuss the proposed particle systems in detail and show their existence (Theorem \ref{thm existence}). As an auxiliary result, we establish a generalization of Schauder's fixed-point theorem (Theorem \ref{thm:schauder}) for mappings defined on the Skorokhod space with the M1 topology.
Moreover, we investigate how the topology of a network affects the occurrence of jumps of a representative particle. Since such jumps may often be interpreted as the events of crises (e.g. a significant number of defaults or infections occurring instantly), we refer to them as the times of fragility and analyze them in detail in Section \ref{se:fragility}. Remarkably, the times of fragility turn out to be closely connected to the Perron-Frobenius eigenvalue of the corresponding graph (see Theorem \ref{cor_PF}).

The second contribution of this work is in the analysis of the associated mean field games. The literature on mean field games is vast (see e.g. \cite{MFG_CarmonaDelarue}, \cite{MFG1}, \cite{MFG2}, and the references therein). However, to date, only \cite{DelarueHAL} has considered mean field games on graphs, using exogenously given homogeneous (random) connection structures. Herein, we allow the particles to control the weights of their connections, thus, obtaining the connection structure endogenously, in equilibrium. In Section \ref{se:game}, we present the general form of the proposed ``resource-sharing" game, which can be viewed as a multi-agent version of the notorious ``network flow'' problem in combinatorial optimization. In such a multi-agent version, the weights of the network edges (i.e. the ``flow") are not determined by optimizing a global criterion, but are obtained from a Nash equilibrium in which every node optimizes its individual criterion locally. In order to obtain concrete results, we focus on a more specific ``credit network" game. There exists a fairly large body of work on modeling equilibrium in credit networks, or, more generally, in asset-liability networks, see e.g. \cite{Atkenson}, \cite{Erol2}, \cite{Acemoglu1}, \cite{Babus}, \cite{Elliott}, \cite{Farboodi}, \cite{Neklyudov}, \cite{Wang}, and the references therein. However, most of the existing models are either static or do not allow for the connections to be determined endogenously in equilibrium\footnote{The exceptions are \cite{Neklyudov}, \cite{Wang}. Nevertheless, the latter models do not include default risk, which is at the center of the present analysis.}. The model proposed herein does possess both of these features. The output of such a model includes the equilibrium dynamics of the interest rate and of the credit connections in the system, which are easy to compute numerically. It turns out that the equilibrium connections are such that the system never exhibits a jump: each particle decreases the strength of its connections to the particles that are close to the barrier, eliminating the cascades. Thus, allowing the particles to act strategically has a regularizing effect on the system. It is worth mentioning that our construction of an equilibrium relies upon the existing results on linear systems in the max-plus algebra (Proposition \ref{prop:game.prop2}), making an interesting connection to a very different area of mathematics.

The rest of the paper is organized as follows. Section \ref{se:existence} introduces the mean field systems with exogenously given connection structures and proves their existence (Theorem \ref{thm existence}). It also contains Theorem \ref{thm:schauder}, which is a generalization of Schauder's fixed-point theorem for mappings defined on the Skorokhod space with the M1 topology. Section \ref{se:fragility} is concerned with the times of fragility of such particle systems. Namely, Theorem \ref{cor_PF} provides necessary and sufficient conditions for fragility in terms of the marginal distributions of the solution, the topology of the network, and the dynamics of the driving processes. In Section \ref{se:game}, we consider a game in which every particle is allowed to control the strength of its connections dynamically, aiming to optimize its objective. The individual optimization problems and the notion of equilibrium are introduced in Subsection \ref{subse:game.setup}. The approach for constructing an equilibrium is described in Subsection \ref{subse:game.approach}. The construction is split into two steps: the dynamic problem (Subsection \ref{subse:game.dynamic}) and the static problem (Subsection \ref{subse:game.static}). The main theorem, showing the existence of an equilibrium and describing its structure, is given in Subsection \ref{subse:game.static} (Theorem \ref{thm:game.main}). The appendix contains technical proofs and auxiliary statements.

\section{Definition and existence of mean field dynamics}
\label{se:existence}

Herein, we deal with a system of particles characterized by their ``level of healthiness" $Y^x$, taking values in $[0,\infty)$, and by their ``type" $x\in{\mathcal X}$, with an abstract finite set $\mX$. We assume that the type of each particle remains constant, whereas the healthiness level changes stochastically over time. When the healthiness level drops to zero, the particle ``dies". In addition, we consider a directed weighted graph, whose nodes are the elements of $\mX$. The connection structure of this graph is encoded by a nonnegative function $C$ and a stochastic kernel $\kappa$, on $\mX$. The weight of the edge from $x$ to $x'$ is given by $C(x)\kappa(x,\{x'\})$. A weight of zero is interpreted as the absence of a connection. The strength of the connection from a particle of type $x$ to the particles of type $x'$ is determined by the weight of the associated edge of the graph. If a particle dies, it shifts down the healthiness levels of all particles that are connected to it. The size of each shift is given in terms of the corresponding connection strength and an interaction function $g$.  In the absence of such shifts, the healthiness level of a particle of type $x$ evolves according to a continuous stochastic process $Z^x$, with $Z^x_0\geq0$.
Thus, $\{Y^x\}_{x\in\mX}$ are defined as a family of stochastic processes satisfying, for each $x\in\mX$,
\begin{equation}\label{what is sol}
\begin{split}
& Y^x_t=Z^x_t+C(x)\int_{\mathcal{X}} g\big(\pp(\tau^{x'}>t)\big)\,\kappa(x,\mathrm{d}x'),\quad t\in[0,\tau^x\wedge T], \\
&Y^x_t=Y^x_{\tau^x},\quad t\in(\tau^x\wedge T,\,T], \\
& \tau^x=\inf\{t\in[0,T]:\,Y^x_t\le 0\}.
\end{split}
\end{equation}
We impose the following conditions on the function $g$ in \eqref{what is sol} throughout.
\begin{assumption}\label{g_ass}
The function $g:\,[0,1]\to[-\infty,0]$ is non-decreasing, on $(0,1]$ it is continuous with values in $(-\infty,0]$, and it is normalized to satisfy $g(1)=0$.
\end{assumption}
 
We point out that the dependence between the driving processes $Z^x$ across $x\in\mathcal{X}$ is immaterial, as each pair $(Y^x,Z^x)$ can be defined on its own probability space.
Note that, although $(\mX,\{C(x)\kappa(x,\{x'\})\}_{x,x'\in\mX})$ is a finite graph, the actual network on which the particles evolve is an infinite (directed, weighted) ``dense" graph. Namely, if $C(x)\kappa(x,\{x'\})>0$, then every particle of type $x$ is connected to all particles of type $x'$ that are still alive (and there is an infinite number of them), with the weight of each connection being the same across all particles of type $x'$.

The interpretation of the proposed system is as follows. The particles may represent cells, individuals, or organizations, having certain ``exposures" to each other. For example, a neuron in a human brain is physically connected to surrounding neurons, a computer in a network allows other computers to upload data to it, and banks may lend funds to other banks and households. The strength of each connection indicates the size of this exposure. When a particle dies, it causes instantaneous shifts in the healthiness levels of the particles connected to it. For example, when the membrane potential in a neuron exceeds a certain threshold, the membrane potentials of the neurons connected to it rise. Once a computer gets infected with a virus, the latter initiates attempts to spread to the connected computers. Similarly, if a company defaults, it causes instantaneous losses to the companies that have lent funds to it.

\medskip

The aforementioned applications to neural and lending networks have been discussed, respectively, in \cite{Delarue1}, \cite{Delarue2}, \cite{Perthame}, \cite{Perthame2}, and in \cite{NadShk}, \cite{HaLeSo}. However, all of the existing models have been limited to fully connected systems, with Brownian dynamics, homogeneous across particles. In the present notation, it corresponds to 
\begin{equation}\label{special case}
\mathcal{X}:=\{0\},\quad Z^0_t:=Z^0_0+\alpha t+\sigma B_t,\;\;t\in[0,T],\quad \kappa(0,\cdot):=\delta_0,
\end{equation}
with $\alpha\in\rr$, $\sigma>0$ and $B$ being a standard Brownian motion. In \cite{NadShk}, $g=\log$ and $Y^0$ stands for the logarithmic monetary reserve process of a representative bank in a large banking system. The profits and losses that the bank makes are due to external investments (giving rise to the $Z^0$ term) and to the defaults of other banks, who have borrowed funds from it (giving rise to the $C(0)\log\pp(\tau^0>t)$ term). For a more detailed explanation, we refer the reader to \cite[Theorem 2.3 and the paragraph preceding it]{NadShk}, as well as the discussion in the first part of Section \ref{se:game}. However, in reality, the borrowing and lending exposures are not distributed uniformly over all members of a credit network. 
For example, the types $x\in\mathcal{X}$ may represent banks' key characteristics such as the geographic region and the size. It is, then, natural to choose the processes $Z^x$, of returns from external investments, as well as the credit exposures, given by $(C(x),\kappa(x,\cdot))$, to depend upon the type. Similarly, the connections between neurons are established according to their physical proximity, or according to the common role played by each cell group. This explains the need for models of the form (\ref{what is sol}), with general $\mX$, $\{Z^x\}_{x\in\mX}$, $\{C(x)\}_{x\in\mX}$, and $\kappa$.

\medskip

The main goal of this section is to establish the existence of a solution to (\ref{what is sol}). Note that, herein, we do not address two important questions: the uniqueness of the solution and the convergence of the associated finite particle systems. The question of uniqueness is a very challenging one. It was only recently resolved by \cite{DelNadShk} in the simpler setting of \eqref{special case}, with a single type and Brownian dynamics. Similarly, the convergence of the finite particle systems has been established only in such a simple setting, see \cite{Delarue2} and \cite{NadShk}. Nevertheless, the methods used to establish the convergence in the latter papers extend to the multi-type setting in a straightforward way, provided that the driving processes $\{Z^x\}_{x\in\mX}$ are Brownian motions. Other types of driving processes may require adjustments in the arguments, or even a brand new strategy of the proof. It is worth mentioning that \cite{NagasawaTanaka1}, \cite{Wu1}, \cite{Wu2} prove the convergence of mean field particle systems whose connection structures are given by a ``dense" graph, as the one proposed herein\footnote{Strictly speaking, \cite{Wu1} studies a sequence of finite particle systems on inhomogeneous random graphs, which converges to a system whose connection structure is similar to the one proposed herein.}. However, their systems are of a classical, smooth, type, without the singular interaction through hitting times. The latter feature makes the problem significantly more complicated, and, at the same time, its presence is necessary to be able to capture certain real world phenomena, as discussed in the next paragraph. 

\medskip

Before we proceed to the existence result, it is important to emphasize the role of the singular interaction through hitting times between the particles. Already in the special case \eqref{special case} with $\alpha\le0$ and $-C(0)\int_0^1 g(z)\,\mathrm{d}z>\E[Z^0_0]$, the paths of the solutions $Y^0$ to \eqref{what is sol} cannot be continuous, as an adaptation of \cite[proof of Theorem 1.1]{HaLeSo} reveals. Since the driving process $Z^0$ is continuous, the discontinuity of $Y^0$ is very surprising at the first glance and, clearly, must be due to the ``self-excitation" term $\pp(\tau^0>t)$ in \eqref{what is sol}. Indeed, as a particle dies, it causes instantaneous shifts in the healthiness levels of the particles connected to it, potentially causing them to die as well. Thus, in a finite particle system (cf. \cite{Delarue2}, \cite{NadShk}), one may observe a cascade, when several particles die at the same time. Even though, as the number of particles grows, the effect of each single particle on the rest of the population does vanish, these cascades do not disappear in the limit, and the resulting solution $Y^0$ does jump (cf. \cite{Delarue2}, \cite{NadShk}). Note that the times of such jumps have a special meaning in the associated applications. These are the times when, in a single instance, the neurons ``synchronize", a macroscopic fraction of computers get infected, or a macroscopic fraction of firms default. As such times are often the subject of practical interest (e.g. as the times of ``crises"), it is important that (i) the model does possess this feature, and (ii) that we are able to analyze it.

\medskip

The above discussion makes it clear that, in the full generality of \eqref{what is sol}, we need to look for solutions $\{Y^x\}_{x\in\mathcal{X}}$ with values in $D([0,T],\rr)^\mathcal{X}$, the product of the spaces of c\`adl\`ag real-valued paths on $[0,T]$. We equip $D([0,T],\rr)^\mathcal{X}$ with the product topology and each factor $D([0,T],\rr)$ with the topology induced by the inclusion
\begin{equation}\label{eq:inclusion}
D([0,T],\rr)\hookrightarrow\hat{D}([0,T+1],\rr),\quad y_t\mapsto\hat{y}_t:=\begin{cases}
y_t\;\;\;\,\text{if}\;\;\;t\in[0,T], \\
y_T\;\;\;\text{if}\;\;\;t\in(T,T+1],
\end{cases}
\end{equation}
where $\hat{D}([0,T+1],\rr)$ is the space of c\`adl\`ag real-valued paths on $[0,T+1]$ left-continuous at $T+1$, endowed with Skorokhod's M1 topology (see any of \cite[Section 4.1]{Delarue2}, \cite{Skorohod}, \cite[Chapter 12]{Whitt} for a detailed discussion of the latter). Our first main result establishes the existence of solutions $\{Y^x\}_{x\in\mathcal{X}}$ to \eqref{what is sol} with values in $D([0,T],\rr)^\mathcal{X}$ under the following assumption.

\begin{assumption}\label{asmp_ex}
\begin{enumerate}[(a)]
\item There exists a family $\{(\eta^x_n)_{n\in\nn}\}_{x\in\mathcal{X}}$ of strictly positive sequences decreasing to $0$ that satisfies
\begin{equation}\label{near0_cond}
\int_{\mathcal X} g\big(\pp\big(Z^{x'}_0>\eta^{x'}_n\big)\!\big)\,\kappa(x,\mathrm{d}x')>-\frac{\eta^x_n}{C(x)},\quad x\in\mathcal{X},\quad n\in\nn.
\end{equation}
\item For all $x\in\mathcal{X}$ and $s\in[0,T]$, the law of $Z^x_s$ possesses no atoms. 
\item For all $x\in\mathcal{X}$, $\varepsilon\in(0,T)$, and $Z^x$-stopping times $\theta^x$, it holds
\begin{equation}
\pp\Big(\min_{s\in[\theta^x,\theta^x+\varepsilon]}\,(Z^x_s-Z^x_{\theta^x})<0\,\Big|\,\theta^x\le T-\varepsilon\Big)=1 
\end{equation}
whenever $\,\pp(\theta^x\le T-\varepsilon)>0$.
\item For all $x\in\mathcal{X}$ and non-increasing $f\in D([0,T],\rr)$ with $f_0=0$, one has 
\begin{equation}
\pp\big(Z^x_s+f_s>0,\;s\in[0,T]\big)>0.
\end{equation}
\end{enumerate}
\end{assumption} 

The condition \eqref{near0_cond} ensures that, initially, not too many particles are about to die, i.e. that the initial distribution of their healthiness does not have too much mass around zero. This assumption is needed in order to ensure that $Y^x$ does not jump at $t=0$. Indeed, it is clear intuitively that the jumps of $Y^x$ occur when too many particles die and drag others with them. The latter, in turn, occurs when too many particles are close to zero (this connection is analyzed rigorously in the next section). And, in order to apply the appropriate fixed-point theorem and to establish the existence of a solution, one needs to ensure that $Y^x$ does not jump at zero (and, hence, $Y^x_{0}=Z^x_0$).
To obtain an even clearer picture, let us consider the special case of \eqref{special case} and $g=\log$. Indeed, inserting both into \eqref{near0_cond} and rearranging the resulting inequality one gets
\begin{equation}\label{complete no casc}
\pp\big(Z_0^0\le\eta^0_n\big)< 1-e^{-\eta^0_n/C(0)},\quad n\in\nn,
\end{equation}
provided $(\eta^0_n)_{n\in\nn}$ is chosen to take values in the set of the continuity points for the cumulative distribution function of $Z^0_0$. At the same time, if the slightly weaker condition (note that $1-e^{-\eta/C(0)}=\eta/C(0)+o(\eta)$ as $\eta\downarrow0$) 
\begin{equation}
\pp\big(Z_0^0\le\eta^0_n\big)< \frac{\eta^0_n}{C(0)},\quad n\in\nn
\end{equation}
is violated, we know from \cite[Proposition 2.5]{NadShk} that a solution $Y^0$ right-continuous at $0$ cannot exist. The condition \eqref{near0_cond} is the network version of \eqref{complete no casc}, for general $g$.

The parts (c) and (d) of Assumption \ref{asmp_ex} are fulfilled, for example, in the continuous semimartingale framework,
\begin{equation}\label{eq:nice}
Z^x_t=Z^x_0+\int_0^t \alpha^x_s\,\mathrm{d}s+\int_0^t \sigma^x_s\,\mathrm{d}B_s,\;\; t\in[0,T],\quad x\in\mathcal{X},
\end{equation}
with $B$ being a Brownian motion (with respect to a filtration supporting $(\alpha^x,\sigma^x)$),
if every $\alpha^x$, $\sigma^x$ lie between two constants (possibly depending on $x$) in $\rr$ and $(0,\infty)$, respectively. This can be deduced by means of the Dambis-Dubins-Schwarz theorem (see e.g. \cite[Chapter 3, Theorem 4.6 and Problem 4.7]{KarShr}) and elementary properties of the standard Brownian motion (see e.g. \cite[Chapter 2, Theorem 9.23]{KarShr} for the part (c) and \cite[Chapter 2, Proposition 8.1 and Theorem 5.12]{KarShr} for the part (d)).

\begin{theorem}\label{thm existence}
Under Assumption \ref{asmp_ex}, the problem \eqref{what is sol} admits a $D([0,T],\rr)^\mathcal{X}$-valued solution $\{Y^x\}_{x\in\mathcal{X}}$, with $\{Y^x_0\}_{x\in\mX}=\{Z^x_0\}_{x\in\mX}$.
\end{theorem}

\begin{rmk}\label{rem:ExistContSol}
Under Assumption \ref{asmp_ex}, for any $t\in(0,T]$, Theorem \ref{thm existence} yields the existence of a solution $\{Y^x\}_{x\in\mathcal{X}}$ on the time interval $[0,t)$.
Assume, further, that $\{Y^x\}_{x\in\mathcal{X}}$ is such that the following time-$t$ version of Assumption \ref{asmp_ex}(a) is satisfied: there exists a family $\{(\eta^x_n)_{n\in\nn}\}_{x\in\mathcal{X}}$ of strictly positive sequences decreasing to $0$ such that
\begin{equation}\label{near0_cond.t}
\int_{\mathcal X} \Big(g\big(\pp\big(A^{x'}\cap\{\chi^{x'}>\eta^{x'}_n\}\big)\!\big)
-g\big(\pp\big(A^{x'}\big)\!\big)\!\Big)\,\kappa(x,\mathrm{d}x')>-\frac{\eta^x_n}{C(x)},\quad x\in\mathcal{X},\quad n\in\nn,
\end{equation}
with $A^x=\{\tau^x\ge t\}$ and $\chi^x=Y^x_{t-}$, $x\in\mX$. Then, the solution $\{Y^x\}_{x\in\mX}$ can be extended from $[0,t)$ to $[0,T]$, so that the resulting solution is continuous at $t$. To see this, we first notice that the proof of Theorem \ref{thm existence} applies, essentially, without changes, to the following modification of \eqref{what is sol}:
\begin{equation}\label{what is sol.A}
\begin{split}
& \widetilde{Y}^x_s=Z^x_s+C(x)\int_{\mathcal{X}} \Big(g\big(\pp\big(A^{x'}\cap\{\widetilde{\tau}^{x'}>s\}\big)\!\big)-g\big(\pp\big(A^{x'}\big)\!\big)\!\Big)\,\kappa(x,\mathrm{d}x'),\quad s\in[0,\widetilde{\tau}^x\wedge T], \\
&\widetilde{Y}^x_s=\widetilde{Y}^x_{\widetilde{\tau}^x},\quad s\in(\widetilde{\tau}^x\wedge T,\,T], \quad
\widetilde{\tau}^x=\inf\{s\in[0,T]:\,\widetilde{Y}^x_s\le 0\},
\end{split}
\end{equation}
yielding the existence of a solution $\{\widetilde{Y}^x\}_{x\in\mX}$, taking values $\{Z^x_0\}_{x\in\mX}$ at $s=0$,  for arbitrary (fixed) events $\{A^x\}_{x\in\mX}$, provided \eqref{near0_cond.t} holds for these $\{A^x\}_{x\in\mX}$ and $\{\chi^x\}_{x\in\mX}=\{Z^x_0\}_{x\in\mX}$, and Assumptions \ref{asmp_ex}(b)--(d) hold.
Next, we choose an arbitrary positive random variable $\widetilde{\chi}$ with an atomless distribution supported away from zero, and apply the latter observation to obtain a solution $\{\widetilde{Y}^x\}_{x\in\mX}$ to \eqref{what is sol.A}, with $Z^x_s$ replaced by\footnote{It is easy to see that Assumptions \ref{asmp_ex}(b)--(d) continue to hold after such a replacement.}
$$
\left(Y^x_{t-}\,\bone_{\{\tau^x\geq t\}} + \widetilde{\chi}\,\bone_{\{\tau^x< t\}}\right) + Z^x_{t+s} - Z^x_t,
$$
and with $A^x=\{\tau^x\geq t\}$,  $x\in\mX$, where each $\tau^x$ is the hitting time associated with $Y^x$. 
Recalling that $\widetilde{Y}^x_0\,\bone_{\{\tau^x\geq t\}}=Y^x_{t-}\,\bone_{\{\tau^x\geq t\}}$, it is easy to see that 
\begin{align*}
Y^x_s\,\bone_{[0,t)}(s) + \big(Y^x_s\,\bone_{\{\tau^x<t\}} + \widetilde{Y}^x_{s-t}\,\bone_{\{\tau^x\geq t\}}\big)\,\bone_{[t,T]}(s),\quad s\in[0,T],\quad x\in\mX
\end{align*}
is an extension of $\{Y^x\}_{x\in\mX}$ to $[0,T]$, continuous at $t$. This observation becomes important in the next section. 
\end{rmk}

\begin{rmk}\label{rem:CommonNoise}
Theorem \ref{thm existence} can be also used to obtain the existence of a solution to the mean field system with common noise associated with \eqref{what is sol}. The latter arises if one introduces a common continuous stochastic process $\overline{Z}$, independent of $\{Z^x\}_{x\in\mathcal{X}}$, and replaces $Z^x$ by $Z^x+\overline{Z}$ and $\pp(\tau^{x'}>t)$ by $\pp(\tau^{x'}>t\,|\,\overline{Z}_s,0\le s\le T)$ in \eqref{what is sol}. Indeed, by conditioning on a path of $\overline{Z}$ the system with common noise reduces to \eqref{what is sol}, with each driving process being the sum of $Z^x$ and the given path of $\overline{Z}$. Then, provided the new driving processes satisfy parts (c) and (d) of Assumption \ref{asmp_ex}, e.g. when $Z^x$, $x\in\mathcal{X}$ and $\overline{Z}$ are ``nice'' semimartingales in the sense of \eqref{eq:nice}, Theorem \ref{thm existence} yields conditional solutions to the mean field system with common noise given almost every path of $\overline{Z}$. Suppose further that there exists a sequence of compact sets $\mathscr{C}_j\subset C([0,T],\rr)$, $j\in\nn$ increasing to a set of full measure for $\overline{Z}$, e.g. by the L\'evy modulus of continuity theorem (see e.g. \cite[Chapter 2, Theorem 9.25]{KarShr}) and the Dambis-Dubins-Schwarz theorem (see e.g. \cite[Chapter 3, Theorem 4.6 and Problem 4.7]{KarShr}) when $\overline{Z}$ is a ``nice'' semimartingale. Then, the set $\mathscr{D}$ in the proof of Theorem \ref{thm existence} can be chosen independently of the path of $\overline{Z}$ on each $\mathscr{C}_j$ and the correspondence sending paths in a given $\mathscr{C}_j$ to the compact range of the respective map $F_2\circ F_1$ therein is weakly measurable by \cite[Theorem 18.5]{Alip} (in our case, the distance function in that theorem is upper semicontinuous in the first argument and continuous in the second argument). Thus, one can apply the measurable maximum theorem (see \cite[Theorem 18.19]{Alip}) to make a measurable selection of the conditional solutions, first on each $\mathscr{C}_j$ and then on a set of full measure for $\overline{Z}$, and compose it with $\overline{Z}$ to find a solution to the mean field system with common noise.
\end{rmk}

The proof of Theorem \ref{thm existence} is given in Appendix A. It is based on a version of the Schauder-Tychonoff fixed-point theorem for the space $D([0,T],\rr)^{\mathcal{X}}$. We emphasize that, even for $\mathcal{X}=\{0\}$, the addition operation in $D([0,T],\rr)^{\mathcal{X}}$ is not continuous (see e.g. \cite[Section 12.7]{Whitt}), so that $D([0,T],\rr)^{\mathcal{X}}$ is not a topological vector space and does not fall under the scope of the classical Schauder-Tychonoff fixed-point theorem in \cite{Tychonoff}. We believe that this auxiliary result, established herein, is interesting in its own right.

\begin{theorem}\label{thm:schauder}
Let $S([0,T],\rr)\subset D([0,T],\rr)$ be the subspace of piecewise constant functions with finitely many discontinuity points and  $\mathscr{D}\subset D([0,T],\rr)^{\mathcal{X}}$ be convex, closed and such that $S([0,T],\rr)^{\mathcal{X}}\cap\mathscr{D}$ is dense in $\mathscr{D}$. Then, every continuous map $F:\,\mathscr{D}\to\mathscr{D}$ whose range has compact closure possesses a fixed-point. 
\end{theorem}

\noindent\textbf{Proof.} The following is an adaptation of the proof of the Schauder-Tychonoff fixed-point theorem for locally convex Hausdorff topological vector spaces in \cite{ParkTan} to the space $D([0,T],\rr)^{\mathcal{X}}$. For  $k\in\nn$, we denote by $U^{1/k}(\mathbf{f})$ the open ball around an $\textbf{f}\in D([0,T],\rr)^{\mathcal{X}}$ of radius $1/k$ in the $\sup$ product metric associated with the strong M1-metric on $D([0,T],\rr)$ (see e.g. \cite[Section 12.3]{Whitt}). The open cover
\begin{equation}
\overline{\mathrm{ran}\,F}\subset\bigcup_{\mathbf{f}\in S([0,T],\rr)^{\mathcal{X}}\cap\mathscr{D}} U^{1/k}(\mathbf{f})
\end{equation} 
has a finite subcover $\overline{\mathrm{ran}\,F}\subset\cup_{j=1}^J U^{1/k}(\mathbf{f}^{j,k})$ due to the assumed compactness of $\overline{\mathrm{ran}\,F}$. Hereby, we have relied on the density of $S([0,T],\rr)^{\mathcal{X}}\cap\mathscr{D}$ in $\mathscr{D}$.  

\medskip

Next, we assign to $\mathbf{f}^{j,k}$, $j=1,2,\ldots,J$ the sets 
\begin{equation}
{\mathscr A}^{j,k}:=\{\mathbf{f}\in \mathscr{D}:\,F(\mathbf{f})\notin U^{1/k}(\mathbf{f}^{j,k})\}\cap\mathrm{conv}\{\mathbf{f}^{1,k},\,\mathbf{f}^{2,k},\,\ldots,\,\mathbf{f}^{J,k}\},\quad j=1,\,2,\,\ldots,\,J,
\end{equation} 
respectively, where $\mathrm{conv}$ stands for the convex hull of a set. We observe that each ${\mathscr A}^{j,k}$ is closed in the Euclidean topology on $\mathrm{conv}\{\mathbf{f}^{1,k},\,\mathbf{f}^{2,k},\,\ldots,\,\mathbf{f}^{J,k}\}$, since a sequence in ${\mathscr A}^{j,k}$ converging in the Euclidean topology 
admits a common finite set of discontinuity points, so that its componentwise convergence in the M1 sense is a consequence of the componentwise pointwise convergence, and the set $\{\mathbf{f}\in \mathscr{D}:\,F(\mathbf{f})\notin U^{1/k}(\mathbf{f}^{j,k})\}$ is closed in the product M1 sense. Moreover, $\cap_{j=1}^J \mathscr{A}^{j,k}=\emptyset$ in view of $\mathrm{ran}\,F\subset\cup_{j=1}^J U^{1/k}(\mathbf{f}^{j,k})$. Thus, the Knaster-Kuratowski-Mazurkiewicz lemma (see \cite{KKM}) demonstrates that, for some $1\le j_1<j_2<\cdots<j_L\le J$ and $\mathbf{f}^{0,k}\in\mathrm{conv}\,\{\mathbf{f}^{j_1,k},\mathbf{f}^{j_2,k},\ldots,\mathbf{f}^{j_L,k}\}$, it holds $\mathbf{f}^{0,k}\notin\cup_{\ell=1}^L \mathscr{A}^{j_\ell,k}$. In other words, $F(\mathbf{f}^{0,k})\in\cap_{\ell=1}^L U^{1/k}(\mathbf{f}^{j_\ell,k})$ and, by the triangle inequality, $\mathbf{f}^{j_{\ell'},k}\in \cap_{\ell=1}^L U^{2/k}(\mathbf{f}^{j_\ell,k})$, $\ell'=1,2,\ldots,L$.

\medskip

Lastly, one can use $\mathbf{f}^{j_\ell,k}\in S([0,T],\rr)^{\mathcal{X}}$, $\ell=1,2,\ldots,L$ and the definition of the strong M1-metric in terms of graph parameterizations to deduce that $U^{2/k}(\mathbf{f}^{j_\ell,k})$, $\ell=1,2,\ldots,L$ are convex, and therefore $\cap_{\ell=1}^L U^{2/k}(\mathbf{f}^{j_\ell,k})$ is convex as well. Putting this together with $\mathbf{f}^{0,k}\in\mathrm{conv}\,\{\mathbf{f}^{j_1,k},\mathbf{f}^{j_2,k},\ldots,\mathbf{f}^{j_L,k}\}$ we end up with $\mathbf{f}^{0,k}\in\cap_{\ell=1}^L U^{2/k}(\mathbf{f}^{j_\ell,k})$, hence, $F(\mathbf{f}^{0,k})\in U^{3/k}(\mathbf{f}^{0,k})$ thanks to the triangle inequality. Taking $k\to\infty$ and appealing to the compactness of $\overline{\mathrm{ran}\,F}$ we infer the existence of a subsequence of $F(\mathbf{f}^{0,k})$, $k\in\nn$ converging to an $\mathbf{f}^{0,\infty}\in\mathscr{D}$. But $\mathbf{f}^{0,k}$, $k\in\nn$ must tend to $\mathbf{f}^{0,\infty}$ along the same subsequence, and the continuity of $F$ yields $\mathbf{f}^{0,\infty}=F(\mathbf{f}^{0,\infty})$. \ep

\medskip

We conclude this section with a discussion of the uniqueness of the solutions to \eqref{what is sol}, which also justifies the definition of the states of fragility given in the next section.
The solutions to \eqref{what is sol} are not unique in general. For instance, consider the situation of \eqref{special case} with $\pp(Z^0_0\in [1,3])=\pp(Z^0_0\in[z_0-1,z_0+1])=1/2$, where $z_0>1-C(0) g(1/2)$ and $-C(0) g(1/2)\ge 3$. Then, Assumption \ref{asmp_ex} is satisfied, as follows from the discussion after the assumption.
Hence, Theorem \ref{thm existence} applies and yields the existence of a $D([0,T],\rr)$-valued solution to \eqref{what is sol}, with $\pp(Y_0^0=Z^0_0>0)=1$. On the other hand, we can obtain a different $D([0,T],\rr)$-valued solution of \eqref{what is sol} as described next. We set $\widetilde{Z}_0^0=Z^0_0+C(0) g(1/2)$, notice that $\PP(\widetilde{Z}^0_0 >0)=1/2$, choose an arbitrary positive random variable $\widetilde{\chi}$ with an atomless distribution supported away from zero, and consider the following equations for $\widetilde{Y}^0$:
\begin{equation*}
\begin{split}
\widetilde{Y}^0_t=\big(\widetilde{Z}_0^0\,\bone_{\{\widetilde{Z}^0_0>0\}} + \widetilde{\chi}\,\bone_{\{\widetilde{Z}^0_0\leq 0\}} \big) +\alpha t+\sigma B_t
+C(0)\big(g\big(\pp(\widetilde{Z}_0^0>0,\,\widetilde{\tau}^0>t\big)\!\big)-g(1/2)\big),\\ 
t\in[0,\widetilde{\tau}^0\wedge T], \\
\widetilde{Y}^0_t=\widetilde{Y}^0_{\widetilde{\tau}^0},\quad t\in(\widetilde{\tau}^0\wedge T,\,T],\quad \widetilde{\tau}^0=\inf\{t\in[0,T]:\,\widetilde{Y}^0_t\le 0\}. \qquad\qquad\qquad\qquad\qquad
\end{split}
\end{equation*}
Since, in the present case, $g(\pp(\widetilde{Z}_0^0>0,\,\widetilde{Y}^0_0>\eta_n^0)\!)-g(\pp(\widetilde{Z}_0^0>0)\!)>-\eta^0_n/C(0)$, for small enough $\eta^0_n>0$, we conclude, as in Remark \ref{rem:ExistContSol}, that the proof of Theorem \ref{thm existence} applies, essentially, without changes and yields the existence of $\widetilde{Y}^0$ solving the above equation, with $\widetilde{Y}^0_0 \,\bone_{\{\widetilde{Z}^0_0>0\}} = \widetilde{Z}_0^0\,\bone_{\{\widetilde{Z}^0_0>0\}}$. Then, it is easy to see that
$$
\widetilde{Z}_0^0\,\bone_{\{\widetilde{Z}^0_0\leq 0\}} + \widetilde{Y}^0_t\,\bone_{\{\widetilde{Z}^0_0>0\}}
$$
is a solution to \eqref{what is sol} with the initial value $\widetilde{Z}^0_0$. Note that the second solution constructed jumps at $t=0$ (more precisely, its time $0$ value does not coincide with $Z^0_0$ which is to be interpreted as the time $0-$ value) even though the distribution of $Z^0_0$ has no mass around zero. Heuristically speaking, a solution may ``jump when it does not have to".
In particular, the jumps of $\{Y^x\}_{x\in\mathcal{X}}$ are not determined intrinsically by \eqref{what is sol}, and one can only hope for the uniqueness of the solution under an extra condition that pins down the jump sizes.
In \cite{Delarue2} and \cite{NadShk}, which study systems of the simpler form \eqref{special case}, the authors choose to define every jump size as follows:
\begin{equation}\label{eq.minJump.1}
Y^0_t=Y^0_{t-}+C(0) \big(g\big(\pp\big(\tau^0\ge t,\,Y^0_{t-}>D_t\big)\!\big)-g\big(\pp(\tau^0\ge t)\!\big)\!\big),\;\;t\in[0,T],
\end{equation}
where
\begin{equation}\label{eq.minJump.2}
\begin{split}
D_t\!:=\!\inf\!\Big\{\!z>0:z+C(0) \big(g\big(\pp\big(\tau^0\ge t,\,Y^0_{t-}>z\big)\!\big)-g\big(\pp(\tau^0\ge t)\!\big)\!\big)>0\Big\},\;\;t\in[0,T].
\end{split}
\end{equation}
This definition corresponds to choosing the \emph{minimal} magnitude of a jump at any time $t$: this is the smallest size of a jump that a solution to \eqref{what is sol} may have at the given time. In the simpler setting of \eqref{special case} and with $g(z)=z-1$ and $g=\log$, respectively, it is shown in \cite{Delarue2} and \cite{NadShk} that such a choice corresponds to a natural definition of a cascade in the associated finite particle systems, and that these particle systems converge to a solution of \eqref{what is sol}, with the jump sizes given by \eqref{eq.minJump.1}--\eqref{eq.minJump.2}.
A direct analogue of the above jump size specification, in the case of general $\mX$, reads:
\begin{equation}
Y^x_t\!=\!Y^x_{t-}+C(x)\!\int_{\mathcal{X}} \!\big(g\big(\pp\big(\tau^{x'}\ge t,\,Y^{x'}_{t-}>D_t\big)\!\big)-g\big(\pp(\tau^{x'}\ge t)\!\big)\!\big)\,\kappa(x,\mathrm{d}x'),\,t\!\in\![0,T],\, x\!\in\!\mathcal{X},
\end{equation}
\begin{equation}
\begin{split}
D_t\!:=\!\inf\!\Big\{\!z\!>\!0:z\!+\!C(x)\!\int_{\mathcal{X}} \big(g\big(\pp\big(\tau^{x'}\!\ge\! t,\,Y^{x'}_{t-}\!>\!z\big)\!\big)\!-\!g\big(\pp(\tau^{x'}\!\ge\! t)\!\big)\!\big)\,\kappa(x,\mathrm{d}x')\!>\!0,\,x\!\in\!\mathcal{X}\Big\}, \\
t\in[0,T].
\end{split}
\end{equation}
However, unlike the case of $\mX=\{0\}$, for a general $\mX$, a heuristic analysis of the associated finite particle systems reveals that other jump size specifications may be relevant. To wit, one may treat the particles of different types asymmetrically when resolving a cascade. For this reason, herein, we do not choose a particular specification of the jump sizes. Instead, for any given $\{Z^x\}_{x\in\mX}$, $\{C(x)\}_{x\in\mX}$, $g$ and $\kappa$, we analyze the times when any solution to \eqref{what is sol} has to jump, which we refer to as the times of fragility (as a jump time is often interpreted as the time of a crisis).


\section{Characterizing the states of fragility}
\label{se:fragility}

In addition to Assumption \ref{g_ass}, we assume throughout this section that, on $(0,1]$, $g$ is continuously differentiable with strictly positive derivative, and that it is normalized to satisfy $g'(1)=1$.\footnote{This normalization causes no loss of generality: it can always be achieved by adjusting $\{C(x)\}_{x\in\mX}$.} Motivated by the discussion at the end of the previous section, herein, we aim to characterize the times of fragility of a given system $(Z,C,g,\kappa)$. The first question we need to address is: in which form do we want to describe these times? Recalling the relevant applications, we notice that it is reasonable to assume that the distribution of every $Y^x_{t-}\,\mathbf{1}_{\{\tau^x\ge t\}}$ is observable at time $t$: indeed, it can be identified with the empirical distribution of the survived particles ``right before" time $t$. In addition, the results of \cite{Delarue2} and \cite{NadShk}, as well as the discussion at the end of the previous section (recall \eqref{eq.minJump.1}--\eqref{eq.minJump.2}), indicate that the jump size of $Y^x$ at time $t$ is determined by the distributions of the survived particles $\mathcal{L}(Y^{x'}_{t-}\,\mathbf{1}_{\{\tau^{x'}\ge t\}})$, across $x'\in\mX$. Thus, for a given $(Z,C,g,\kappa)$, we say that a time $t\in[0,T]$ and a collection of probability distributions $\mathcal{L}(Y^{x'}_{t-}\,\mathbf{1}_{\{\tau^{x'}\ge t\}})$ on $[0,\infty)$ form a fragile state, if any solution to \eqref{what is sol}, with the given marginal distributions at time $t$, has a non-zero jump at that time. In fact, it is easy to guess that only the asymptotic properties of the normalized marginal distributions near zero determine whether a jump occurs or not. Hence, for any $\mathcal{L}(Y^{x}_{t-}\,\mathbf{1}_{\{\tau^{x}\ge t\}})$, we introduce the associated $\{c^x\geq0\}_{x\in\mathcal{X}}$ and $\{z^x>0\}_{x\in\mathcal{X}}$, such that
\begin{equation}\label{eq.cz.def}
\pp\big(\tau^x\ge t,\,Y^x_{t-}\in(0,z)\!\big)\,g'\big(\pp(\tau^x\ge t)\!\big)\ge c^x z,\;\; z\in[0,z^x],\quad x\in\mathcal{X}.
\end{equation}
Clearly, such $\{c^x\}_{x\in\mX}$, $\{z^x\}_{x\in\mX}$ always exist, and there are typically infinitely many choices of such constants. The larger $c^x$ and $z^x$ are, the more particles of type $x$ are concentrated around zero (i.e. the more fragile the state is). Indeed, if $\mathcal{L}(Y^x_{t-}\,\mathbf{1}_{\{\tau^x\ge t\}})$ has a density $p^x$ continuous at zero, then $c^x$ can be any number in $(0,p^x(0)g'(\pp(\tau^x\ge t)\!))$, and the closer $c^x$ is to $p^x(0)g'(\pp(\tau^x\ge t)\!)$, the smaller is $z^x$. To develop a better intuition for the statements that follow, it is convenient to interpret $c^x$ as $p^x(0)g'(\pp(\tau^x\ge t)\!)$: indeed, $c^x$ is simply a proxy for $p^x(0)g'(\pp(\tau^x\ge t)\!)$ in the cases where the latter is not well-defined.

\smallskip

In this section, for a general quadruple $(Z,C,g,\kappa)$, we formulate a sufficient condition for a state to be fragile, expressed as the existence of associated families $\{c^x,z^x\}_{x\in\mX}$ satisfying additional properties. We show how this sufficient condition can be used and that it is fairly sharp.
As a preparation for the upcoming Theorem \ref{cor_PF}, we supply the following assumption. 

\begin{assumption}\label{asmp_dom}
For any $t\in[0,T)$ and any $\{z^x>0\}_{x\in\mathcal{X}}$, there exist a constant $\theta\in(0,T-t)$ and a random field $\xi^x_{t_1,t_2}$, $t\le t_1<t_2\le t+\theta$, $x\in\mathcal{X}$, such that:
\begin{enumerate} 
\item[(a)] for every $z\in(-\infty,z^x]$, $t\le t_1<t_2\le t+\theta$, and $x\in\mathcal{X}$,
\begin{equation}\label{eq:dom}
\pp\Big(-\min_{t_1\le s\le t_2} (Z^x_s-Z^x_t)\ge z\,\Big|\,Z^x_s,\,s\in[0,t]\Big)
\ge\pp(\xi^x_{t_1,t_2}\ge z),
\end{equation}
\item[(b)] and, for at least one sequence $s_n\downarrow t$,
\begin{equation}
\lim_{n\to\infty}\,\frac{1}{n}\log\big(\min_{x\in\mathcal{X}}\, \E[\xi^x_{s_{n+1},s_n}\wedge z^x]\big)\ge0. \label{regular2}
\end{equation}
\end{enumerate}
\end{assumption}

The first part of the above assumption states that the absolute size of the drawdowns of $Z^x$ can be dominated stochastically from above, uniformly over the paths of $Z^x$ up to time $t$. The second part states that these drawdowns do not decay too fast as the time horizon decreases, so that the particles drop sufficiently below the initial level $Y^x_{t-}\,\mathbf{1}_{\{\tau^x\ge t\}}$. It is shown in the proof of Theorem \ref{cor_semimart}, in Appendix B, that Assumption \ref{asmp_dom} is satisfied by any $Z$ in the form
\begin{equation}\label{semimart_set}
Z^x_s=Z_0^x+\int_0^s \alpha^x_r\,\mathrm{d}r+\int_0^s \sigma^x_r\,\mathrm{d}B_r,\;\;s\in[0,T],\quad x\in\mathcal{X},
\end{equation} 
with $\alpha^x\le\overline{\alpha}$ and $\underline{\sigma}\le \sigma^x\le\overline{\sigma}$, for all $x\in\mathcal{X}$, with suitable $\overline{\alpha}\in\rr$ and $\underline{\sigma},\overline{\sigma}\in(0,\infty)$.
On the other hand, outside of the class of semimartingales with bounded characteristics, there exist processes that do not satisfy Assumption \ref{asmp_dom}. An example is presented at the end of this section.

\begin{definition}
Let the matrix $(C(x)\kappa(x,\{x'\})c^{x'})_{x,x'\in\mathcal{X}}$ be irreducible. The limit
\begin{equation}\label{what is rho}
\varrho:=\lim_{n\to\infty}\, \frac{1}{n}\log\int_{\mathcal{X}^n} \prod_{m=0}^{n-1} (C(x_m)c^{x_{m+1}})\,\kappa(x_0,\mathrm{d}x_1)\,\kappa(x_1,\mathrm{d}x_2)\,\ldots\,\kappa(x_{n-1},\mathrm{d}x_n)
\end{equation} 
is referred to as its \emph{logarithmic Perron-Frobenius eigenvalue} (see the Perron-Frobenius theorem in the form of \cite[Theorem 3.1.1(e)]{DZ}).
\end{definition}

To understand the purpose of the above definition, we need to recall the existing results on fragility of single-type systems, developed in \cite{Delarue2}, \cite{NadShk}.
Namely, assuming $\mX=\{0\}$ and $g(z)=z-1$ or $g=\log$, it is shown in the latter papers that a sufficient condition for fragility of a given state is the existence of $(c^0,z^0)$, satisfying (\ref{eq.cz.def}), such that $C(0) c^{0}\geq1$.\footnote{In fact, the affine choice of $g(z)=z-1$ in the model of \cite{Delarue2}  allows the authors to show that this condition is both necessary and sufficient.}
The interpretation of this condition is also clear. Recall that $c^0$ is a proxy for the (normalized) density of the particles that are about to hit zero. If this density reaches the critical level $1/C(0)$, which is inversely proportional to the connectivity level $C(0)$ among the particles, a cascade occurs. It is also easy to notice that, in the single-type case, $\varrho = \log(C(0) c^{0})$. Thus, the sign of $\varrho$ determines the fragility of a given state. In a multi-type setting, $C(x)c^{x}$ may exceed $1$ for some $x\in\mX$, and it may be below $1$ for other $x$. The Perron-Frobenius eigenvalue combines the ``loss intensities" $\{C(x)c^x\}_{x\in\mathcal{X}}$ with the topology of the network, represented by $\kappa$, to determine whether a given state is fragile.

\begin{theorem}\label{cor_PF}
Fix an arbitrary $t\in[0,T)$ and assume that $\{Z^x\}_{x\in\mathcal{X}}$ satisfy Assumptions \ref{asmp_ex}(b)--(d) and \ref{asmp_dom}.
\begin{enumerate}
\item[(a)] Suppose that a solution $\{Y^x\}_{x\in\mathcal{X}}$ of \eqref{what is sol} admits $\{c^x, z^x\}_{x\in\mX}$ satisfying \eqref{eq.cz.def} and is such that at least one closed irreducible component of the matrix $(C(x)\kappa(x,\{x'\})c^{x'})_{x,x'\in\mathcal{X}}$ has a strictly positive logarithmic Perron-Frobenius eigenvalue $\varrho$. Then, $Y^x_{t-}>Y^x_t$ for at least one $x\!\in\!\mathcal{X}$.
\item[(b)] Suppose that a solution $\{Y^x\}_{x\in\mathcal{X}}$ of \eqref{what is sol} admits $\{c^x, z^x\}_{x\in\mX}$ satisfying  
\begin{equation}
\pp\big(\tau^x\ge t,\,Y^x_{t-}\in(0,z)\!\big)\,g'\big(\pp(\tau^x\ge t)\!\big)\le c^x z,\;\; z\in[0,z^x],\quad x\in\mathcal{X},
\end{equation}
and is such that every closed irreducible component of the matrix $(C(x)\kappa(x,\{x'\})c^{x'})_{x,x'\in\mathcal{X}}$ has a strictly negative logarithmic Perron-Frobenius eigenvalue $\varrho$.
Then, there exists a solution $\{\tilde{Y}^x\}_{x\in\mathcal{X}}$ of \eqref{what is sol} coinciding with $\{Y^x\}_{x\in\mathcal{X}}$ on $[0,t)$ such that $\tilde{Y}^x_{t-}=\tilde{Y}^x_t$ for every $x\in\mathcal{X}$.
\end{enumerate}
\end{theorem}

\begin{rmk}
If the distributions $\{\mathcal{L}(Y^x_{t-}\,\mathbf{1}_{\{\tau^x\ge t\}})\}_{x\in\mX}$ have densities $\{p^x\}_{x\in\mX}$ continuous at zero, then the above theorem implies that: (a) this state is fragile if $\varrho>0$ for at least one closed irreducible component of $(C(x)\kappa(x,\{x'\})p^{x'}(0)g'(\pp(\tau^{x'}\ge t)\!))_{x,x'\in\mathcal{X}}$, and (b) this state is not fragile if $\varrho<0$ for every closed irreducible component of $(C(x)\kappa(x,\{x'\})p^{x'}(0)g'(\pp(\tau^{x'}\ge t)\!))_{x,x'\in\mathcal{X}}$.
\end{rmk}

\begin{rmk}
The case $\varrho=0$ is more subtle. Theorem \ref{cor_semimart}, in Appendix B, provides a sufficient condition for fragility that covers this case. However, this condition is less explicit, hence, we postpone it to the appendix.
\end{rmk}

The proof of Theorem \ref{cor_PF} is stated in Appendix B. It is worth mentioning that the proof is based on the key Lemma \ref{thm:jump_crit}, which is a generalization of \cite[Proposition 2.7]{Delarue2}.


\smallskip

Theorem \ref{cor_PF} allows us to ignore the exact dynamics of $Z$ (provided it is a ``nice" continuous semimartingale) and to describe the states of fragility only in terms of the current marginal distributions of the solution and the topology of the network. This leads to very interesting observations, such as the following example.
Consider two particle systems, $(Z,C,g,\kappa^1)$ and $(Z,C,g,\kappa^2)$, with arbitrary $\{C(x)>0\}_{x\in\mX}$ and $g$ (the same for both models), with $\kappa^1(x,\mathrm{d}x')=\delta_x(\mathrm{d}x')$, and with $\kappa^2(x,\mathrm{d}x')=\mathrm{d}x'/|\mX|$ (on $\mX$, we denote by ``$\mathrm{d}x'$'' the counting measure). We refer to the first system as the ``completely clustered" network, and to the second one as the ``fully connected" one. The rationale behind these terms is clear: the particles in the first system are only connected within the same type, while the particles in the second system are connected to every other particle in the network, uniformly. Note that the total strength of the out-going connections from each particle is the same in both networks -- only the distributions of the connections are different. Let us also assume that $Z^x$, $x\in\mX$ are the same for both systems and are ``nice" semimartingales, so that Assumption \ref{asmp_dom} is satisfied, and that $\{Z^{x}_0\}_{x\in\mX}$ have densities $\{p^x\}_{x\in\mX}$, continuous at zero.
Thus, with everything else being equal, we can compare the fragility of the two systems at time $t=0$ by computing the associated Perron-Frobenius eigenvalues.
For the first network, the associated matrix $(C(x)\kappa^1(x,\{x'\})p^{x'}(0))_{x,x'\in\mathcal{X}}$ is, clearly, reducible, hence, we analyze every closed irreducible component $\{x\}$ separately\footnote{Strictly speaking, $\{x\}$ may not be irreducible if $C^xp^{x}(0)=0$, but it is easy to see that our conclusion about fragility remains valid even when such cases are taken into account.}:
\begin{equation*}
\varrho^1= \lim_{n\to\infty}\, \frac{1}{n}\log (C(x)p^{x}(0))^n = \log (C(x) p^x(0)).
\end{equation*} 
We conclude that the first system is in a fragile state if $\max_{x\in\mX} C(x) p^x(0)>1$, and the state is not fragile if $\max_{x\in\mX} C(x) p^x(0)<1$.
For any closed irreducible component $\tilde{\mX}$ of the second system, we have:
\begin{equation*}
\varrho^2\!=\!\lim_{n\to\infty}\, \frac{1}{n}\log\int_{\tilde{\mathcal{X}}^n} \prod_{m=0}^{n-1} (C(x_m)p^{x_{m+1}}(0))\,\kappa(\mathrm{d}x_1)\,\ldots\,\kappa(\mathrm{d}x_n)
\!=\!\log \bigg(\frac{1}{|\mX|} \sum_{x\in\tilde{\mathcal{X}}} C(x)p^{x}(0)\!\bigg).
\end{equation*} 
Thus, the second system is in a fragile state if $\max_{\tilde{\mX}}\frac{1}{|\mX|} \sum_{x\in\tilde{\mathcal{X}}} C(x)p^{x}(0)>1$, and the state is not fragile if $\max_{\tilde{\mX}}\frac{1}{|\mX|} \sum_{x\in\tilde{\mathcal{X}}} C(x)p^{x}(0)<1$, where the maximum is taken over all closed irreducible components $\tilde{\mX}$.
Since the average never exceeds the maximum, the above analysis indicates that the fully connected system is more stable (up to the case $\varrho^1=\varrho^2=0$).
Using the existence result in Theorem \ref{thm existence}, one can easily show that the fully connected network is strictly more stable, i.e. there exist states at which the first system is fragile, while the second one is not.


\medskip

We conclude this section by showing the sharpness of Assumption \ref{asmp_dom}. Namely, we present a (fairly simple) stochastic process $Z$ that does not satisfy this assumption and for which Theorem \ref{cor_PF} does fail.
Consider $\mathcal{X}:=\{0\}$ and $Z_t^0:=Z^0_0+f_t+B_t$, $t\in[0,T]$ with a continuous function $f:\,[0,T]\to[0,\infty)$ satisfying $\theta:=\liminf_{t\downarrow 0} (f_t/\sqrt{t})>0$.
It is straightforward to construct $\xi^x_{t_1,t_2}$ satisfying \eqref{eq:dom}, as $Z^0$ has independent increments. Nevertheless, we claim that, regardless of the exact choice of $\xi^x_{t_1,t_2}$, \eqref{regular2} does not hold at $t=0$ as long as $\xi^x_{t_1,t_2}$ satisfies \eqref{eq:dom}. Indeed, for any $z^0>0$, any $(s_n)_{n\in\nn}$ decreasing to $0$, and $\theta_n:=\inf_{t\in(0,s_n]} (f_t/\sqrt{t})\ge0$, $n\in\nn$, the estimates
\begin{equation}
\begin{split}
\E[\xi^0_{s_{n+1},s_n}\wedge z^0] 
& \le \E\big[\!-\!\min_{s_{n+1}\le s\le s_n} (f_s+B_s)\big] \\
& \le -\theta_n\sqrt{s_n}+\E\big[\!-\!\min_{s_{n+1}\le s\le s_n} (\theta_n\sqrt{s_n-s}-B_{s_n}+B_s)\big] \\
& = -\theta_n\sqrt{s_n}+\sqrt{s_n-s_{n+1}}\;\E\big[\!-\!\min_{t\in[0,1]} (\theta_n\sqrt{t}+B_t)\big]
\end{split}
\end{equation}
reveal that either $\E[\xi^0_{s_{n+1},s_n}\wedge z^0]\le 0$ infinitely often or
\begin{equation}\label{geom_series}
s_{n+1}< s_n\bigg(1-\frac{\theta_n^2}{\E[-\min_{t\in[0,1]} (\theta\sqrt{t}/2+B_t)]^2}\bigg) 
\end{equation}
when $n\in\nn$ is sufficiently large. In the first case, it is immediate that \eqref{regular2} is violated, whereas in the second case this can be seen from the fact that $(s_n)_{n\in\nn}$ decays exponentially fast and  
\begin{equation}
\E[\xi^0_{s_{n+1},s_n}\wedge z^0] 
\le\E\big[\!-\!\min_{s_{n+1}\le s\le s_n} \big(\theta_n\sqrt{s}+B_s\big)\big]
\le\sqrt{s_n}\,\E\big[\!-\!\min_{t\in[0,1]} \big(\theta\sqrt{t}/2+B_t\big)\big].
\end{equation}
Let us show that, with such a choice of $Z$, even if $\varrho>0$, the state may not be fragile.
Indeed, consider $f_t=-C(0) g(\pp(Z_0^0+\min_{s\in[0,t]} B_s>0))$, $t\in[0,T]$, and assume that $Z^0_0$ is independent of $B$ and satisfies $\pp(Z_0^0\in(0,z))\ge h z$, $z\in[0,\eps]$ for some $h>1/C(0)$ and $\eps>0$. It is easy to see that $\liminf_{t\downarrow 0} (f_t/\sqrt{t})>0$. 
Choosing $t=0$, $c^0=h$, and $z^0=\eps$, we notice that \eqref{eq.cz.def} is satisfied and that $\varrho=\log (C(0) c^0) >0$. Nevertheless, \eqref{what is sol} does have a continuous solution
\begin{equation}
Y^0_t:=Z_0^0+B_t=Z_0^0+f_t+B_t+C(0)g\big(\pp\big(Z_0^0+\min_{s\in[0,t]} B_s>0\big)\!\big),\quad t\in[0,T].
\end{equation}
Intuitively, the solution $Y^0$ does not jump at $t=0$ because the driving process $Z^0$ has an ``infinitely" large positive drift at that time, which pushes the particles up very fast. Hence, even though a sufficiently large fraction of the particles is concentrated around zero (which would be enough to generate a jump if $Z^0$ was a ``nice" continuous semimartingale), it does not cause a non-zero fraction of them to die simultaneously.

\section{Controlled mean field dynamics and equilibrium}
\label{se:game}

\subsection{Model setup}
\label{subse:game.setup}

In this section, we consider particle systems of the form \eqref{what is sol}, in which the strength of connections across types, given by the kernel $\kappa(x,\mathrm{d}x')$, is controlled dynamically by the particles. We choose natural objectives for the individual particles and analyze the equilibrium connections that arise in the associated mean field game. Remarkably, we show that there exists an equilibrium in which the particles manage to adjust their connections in such a way that a cascade never occurs.

As in the previous sections, we consider an abstract (non-empty) finite set $\mX$. We endow $\mX$ with a probability measure $\mu$ such that $\mu(\{x\})>0$ for all $x\in\mX$. We fix an arbitrary function $O:\,\mX^2\rightarrow\{0,1\}$, which indicates the types that are ``physically" connected to each other, i.e. $x$ is physically (directly) connected to $y$ if and only if $O(x,y)=1$.\footnote{In what follows, we drop the adverb ``directly", and say that two nodes are connected if they are connected by a single edge of the graph (unless stated otherwise).} The pair $(\mX,O)$ defines a directed graph, whose edges indicate the physical connections.
As before, we consider a system of particles characterized by their type $x\in\mX$ and their level of healthiness $Y^x$ (the latter changes stochastically over time).
We assume that there exist infinitely many particles of each type. The value of $\mu(\{x\})$ represents the total mass of all particles of type $x$, and the initial distribution of particles across their characteristics, $(x,y)$, is given by $\PP(Y^x\in\mathrm{d}y)\,\mu(\mathrm{d}x)$.

In what follows, we consider a resource-sharing dynamic network game with a continuum of players. The proposed game can be viewed as a dynamic multi-agent version of the network flow problem (cf. \cite{NetworkFlow}). In this version, instead of choosing a static network topology that optimizes a given global objective, we search for the topology that arises in equilibrium, when each node optimizes its own objective dynamically.
More specifically, the proposed model falls within the class of mean field games (see Section \ref{se:intro} for more on such games).
In this game, the edges of $(\mX,O)$ indicate the opportunities for sharing a given resource, i.e. if $x$ is connected to $x'$, then a particle of type $x$ is allowed to give access to a certain resource for the particles of type $x'$.
In addition, the edges of $(\mX,O)$ have weights, which indicate the amount of resource shared between the particles of two types. The particles are allowed to change these weights dynamically, in order to optimize their objective.
In such a game, the dynamics of a representative particle are given by
\begin{equation}\label{eq.game.Yx.def}
Y^x_t = Z^x_t + \int_{\mX} \int_0^t C^l_{s} \nu_{s}(x') \,\mathrm{d} g(\theta_{s}(x')\!)\, \mu(\mathrm{d}x'),\quad t\geq0,\,\,\,\,x\in\mX,
\end{equation}
where
\begin{equation}\label{eq.game.Zx.def}
Z^x_t := Y^x_0 + \sigma\, W_t
+ \int_0^t L\left[x,\,P(x,C^b_s,C^l_s)
-  C^b_s r_s(x)
+ \int_{\mathcal{X}} C^l_s  r_s(x')\nu_s(x')\mu(\mathrm{d}x')\right] \,\mathrm{d}s,
\end{equation}
$Y^x_0$ is independent of the Brownian motion $W$, $g$ satisfies Assumption \ref{g_ass}, and the following consistency requirements hold:
\begin{equation}\label{eq.game.theta.1}
\begin{split}
& \tau^x=\inf\{t\geq0:\,Y^x_t\leq0\},
\quad\theta_t(x) := \PP\left(\tau^x>t\right) = \PP\Big(\inf_{s\in[0,t]}Y^x_s>0\Big), \\
& O(x,x')=0\,\Rightarrow\,\nu_\cdot(x')=0,\;\; x'\in\mX.
\end{split}
\end{equation}
If the input elements $(C^l=C^l(x),\nu(\cdot)=\nu(x,\cdot))$ are deterministic, constant in time, and the same within each type $x$ (but, possibly, different across $x$), the equations (\ref{eq.game.Yx.def}) and (\ref{eq.game.theta.1}) coincide with \eqref{what is sol}.
The main difference in this part of the paper is that the quantities $(C^b,C^l,\nu(\cdot),r(\cdot))$ are determined endogenously.

Let us discuss the general interpretation of the above dynamics.
Recall that each particle can decide on how much of the resource it is willing to share with other types, thus, creating the weights of the graph edges\footnote{If a particle shares a unit of resource with a certain type, it is distributed uniformly across all particles of this type.}. The number $C^l_t$ denotes the overall amount of resource that the given particle is willing to share at time $t$.\footnote{More precisely, as we consider a continuum player game (cf. \cite{Aumann}), $C^l_t$ denotes the amount of resource per unit mass of agents lent to by the particle. The same applies to other amounts (of resource or currency) related to an individual particle.} The probability measure $\kappa_t(\mathrm{d}x')=\nu_t(x') \mu(\mathrm{d}x')$ denotes the distribution of the total amount $C^l_t$ among various neighbors\footnote{The decomposition $\kappa_t=\nu_t\mu$ makes it more convenient to formulate the market clearing condition (\ref{eq.game.clearingCond}).}, i.e. the contribution of this particle to the weight of the edge from $x$ to $x'$ is $C^l_t \nu_t(x')\mu(\mathrm{d}x')$. Similarly, $C^b_t$ denotes the overall amount of resource it is willing to receive at time $t$.\footnote{Of course, $C^b$ cannot be chosen arbitrarily - it depends on $C^l$ and $\nu$ of other particles. In our definition of equilibrium, we add such a ``clearing" condition.} Gaining access to the resource is costly, i.e. a particle of type $x'$ has to pay $r(x')$ units of currency for having access to one unit of resource per unit of time. Similarly, the particles sharing the resource with a particle of type $x'$ receive $r(x')$. This explains the last two terms on the right-hand side of (\ref{eq.game.Zx.def}). Having received $C^b_t$ and shared $C^l_t$ units of resource at time $t$ yields the profit or loss of $P(x,C^b_t,C^l_t)$ per unit of time, through an internal mechanism of each particle (e.g. the business profits or losses)\footnote{In the simplest case, $P(x,C^b_t,C^l_t)$ is proportional to $C^b_t-C^l_t$, but one can easily come up with other relevant examples where the dependence is non-linear.}.
We assume that the profits and losses affect the healthiness of a particle through the function $L$ (e.g., a company may spend a fraction of its profits on improving its healthiness level), which explains the third term on the right-hand side of (\ref{eq.game.Zx.def}).
The Brownian noise in (\ref{eq.game.Zx.def}) represents the changes in the healthiness of a particle that are due to external factors.
Finally, the death of each particle (occurring when $Y$ hits zero) causes immediate damage to the healthiness of the particles that share the resource with it. The magnitude of this damage is determined by the strength of the connection and by the choice of function $g$, which explains the last term on the right-hand side of (\ref{eq.game.Yx.def}) (see more on the interpretation of this term in the discussion following equation \eqref{what is sol}).
When a particle dies, it leaves the system, in the sense that it can no longer be connected to other particles. The value of $\theta_t(x)$ indicates the fraction of particles of type $x$ that have not died by time $t$.
The processes $(C^b,C^l,\nu)$ are controlled by the particle.

There are several concrete examples where such games arise naturally. For example, consider a cybersecurity problem, where the particles represent companies. Each company may allow others to gain access to its computer resources (e.g. servers, advertisement platforms, and so on), collecting a fee for this service. The other companies agree to pay the fee as this service may generate additional profits for their business (e.g. increase the sales volume). The process $Y$, in this case, represents the level of protection of a given company against cyber-attacks. When $Y$ falls to zero, the company gets infected with a malicious software, which spreads to other companies connected to it, instantaneously decreasing their level of protection. Clearly, there exist many other examples.
However, in order to obtain stronger results, we need to choose a specific interpretation of such a game, leading to a specific choice of the functions $L$ and $P$. In this section, we focus on the credit network game, in which
\begin{equation*}
L(x,z)=z,\quad P(x,c^b,c^l) = \alpha'(x) + \alpha(x)(c^b - c^l),
\quad g=\log,
\end{equation*}
with some $\alpha',\alpha\in[0,\infty)^{|\mX|}$, fixed throughout the rest of the section.
In this game, the particles represent individuals or organizations engaged in lending and borrowing of funds, and $Y^x$ represents the logarithm of the capital\footnote{Strictly speaking, we normalize the capital of each particle so that its default barrier is at $1$.} of a representative particle of type $x$. Whenever $Y^x$ drops to zero, the corresponding particle defaults (i.e. ``dies"), causing immediate losses to its neighbors. Note that the choice $g=\log$ causes no loss of generality: indeed, any $g\in C^1((0,1])$ with $g'>0$ would yield the same right-hand side in \eqref{eq.game.Yx.def}, provided $(C^l,\kappa)$ are adjusted accordingly (recall that the latter are controlled dynamically by the particles). The specific choice $g=\log$ is motivated by the model of \cite{NadShk}, which shows that it is natural when $C^l$ is constant in time and $\mathcal{X}$ is a singleton.
We assume that a regulator has imposed a constraint on the maximum size of non-core liabilities, i.e. the maximum size of $C^l-C^b$, which represents each particle's own funds lent within the system, measured as a multiple of its capital, cannot exceed a given constant $\bC\in(0,\infty)$. Note that such a constraint, effectively, limits the strengths of connections within the system. Our main goal is to determine the connections and the interests rates arising in equilibrium.\footnote{One may also add a constraint on the maximum value of $C^b$, accompanied with a rule (endogenous or exogenous) for how the limited lending opportunities are shared among multiple particles lending to the same type. With such an extension, the particles could benefit from the deaths of others, as the latter may generate additional lending opportunities. We thank the anonymous referee for this observation, but we do not pursue this direction herein, as it would complicate the notation and derivations without providing significant new insights into the main question of our analysis.}


Let us introduce the mathematical conventions we use throughout the rest of the section. We let $\RR_+:=[0,\infty)$ and $\RR_-:=(-\infty,0]$. For any $p\in[0,\infty]$ and any subsets $A$ and $B$ of the Euclidean space, we denote by $L^p(A\rightarrow B)$ the space of all measurable functions, satisfying the associated integrability or essential boundedness conditions. By $\mathbb{L}^{p}(A\rightarrow B)$ we understand the space of equivalence classes of the elements of $L^p(A\rightarrow B)$. We also use the short-hand notations, $L^p(A)$ and $\mathbb{L}^p(A)$, if the associated functions take values in $\RR$, as well as $L^p_+(A)$, $\mathbb{L}^p_+(A)$, $L^p_-(A)$, and $\mathbb{L}^p_-(A)$, if the associated functions take values in $\RR_+$ and $\RR_-$, respectively. Finally, we adopt the convention that all suprema and maxima over the empty set are equal to $-\infty$.

We assume that $C^b,C^l\in L^{0}_+(\RR_+)$; $\nu\in (L^{0}_+(\RR_+))^{|\mX|}$ is such that $\int_{\mX}\nu_t(x)\mu(\mathrm{d}x)\in\{0,1\}$, for all $t\geq0$; $r\in (L^{0}_+(\RR_+))^{|\mX|}$; and $\theta\in (L^{0}_+(\RR_+))^{|\mX|}$ is such that, for all $x\in\mX$, $\theta_{\cdot}(x)$ is a non-increasing function with values in $(0,1]$, starting from $1$.
For convenience, we introduce special notations for the spaces of possible $\theta$ and $\nu$:
\begin{equation*}
\begin{split}
& \Theta = \left\{ \theta\in L^{0}\left(\RR_+\rightarrow(0,1]^{|\mX|}\right):\, \theta_{\cdot}(x)\text{ is non-increasing and }\theta_0(x)=1,\text{ for all $x\in\mX$}\right\}, \\
& \mV_0 = \left\{ \nu\in \RR^{|\mX|}_+ :
\,\int_{\mX}\nu(x')\mu(dx')=1\right\}, \\
& \mV = \left\{ \nu\in L^{0}\big(\RR_+\rightarrow\RR^{|\mX|}_+\big) :
\,\nu_t\in \mV_0,\text{ for almost every }t\geq0\right\}.
\end{split}
\end{equation*}
\begin{rmk}
Note that we have restricted $\theta_t(x)$ to be strictly positive. We could easily lift this restriction and allow for $\theta_t(x)=0$ by making the notation slightly more complicated. We choose not to do it, because in the equilibrium we construct $\theta_t(x)$ does remain strictly positive at all times.

Similarly, we could allow the controls of the particles to be random (e.g. depending on $Y$). However, even then, there would still exist an equilibrium with deterministic optimal controls. Hence, we a priori restrict the setting to deterministic controls, to simplify the notation.
\end{rmk}

\begin{definition}\label{def:game.admis}
For a given pair $(r,\theta)\in (L^{0}_+(\RR_+))^{|\mX|}\times\Theta$, the control $(C^b,C^l,\nu)\in (L^{0}_+(\RR_+))^2\times \mV$ is admissible for a particle of type $x$, if
\begin{enumerate}
\item the integrals on the right-hand sides of (\ref{eq.game.Yx.def}) and (\ref{eq.game.Zx.def}) are finite for every $t\geq0$;
\item $C^l_t \leq C^b_t + \bC$, for almost every $t\geq0$;
\item for almost every $t\geq0$ and all $y\in\mX$, $O(x,y)=0$ implies $\nu_t(y)=0$.
\end{enumerate}
\end{definition}

The second condition in the above definition represents the constraint on the non-core exposure. The third condition ensures that the lending and borrowing occurs only between the particles that are physically connected (according to the adjacency matrix $O$).

\begin{rmk}\label{rem:game.CompetitiveBorrowers}
Note that there is a certain asymmetry between lending and borrowing reflected in the dynamics (\ref{eq.game.Yx.def})--(\ref{eq.game.Zx.def}) and in the definition of admissibility of the controls. Namely, we assume that the particles control both the size and the distribution of their lending, but they only control the size of their borrowing.
Nevertheless, in equilibrium, the lending and borrowing distributions will be matched, so that the market ``clears".
The asymmetry between lending and borrowing is an intentional modeling choice we make, which persists throughout the remainder of the section. Intuitively, this choice corresponds to ``dominant lenders" (equivalently, ``competitive borrowers"), i.e. in any borrowing-lending deal, the lender decides whom to lend (more precisely, to which type to lend), and the borrower accepts the funds from whoever is willing to lend to it. We believe that this is a natural assumption, and it allows us to narrow down the set of possible equilibria, as will be shown later in the paper.
\end{rmk}

To emphasize the dependence of the state process $Y$ on the market and on the controls, we will often write $Y^{x}(r,\theta;C^b,C^l,\nu)$.
A particle with (initial) characteristics $(x,y)$ aims to maximize the following objective:
\begin{equation}\label{eq.game.J.def}
J(x,y;r,\theta;C^b,C^l,\nu):=\EE\bigg[ \int_{0}^{\tau^x} e^{-\gamma t}\,Y^x_t(r,\theta;C^b,C^l,\nu)\,dt \,\bigg|\, Y^x_0=y\bigg],
\end{equation}
with a constant $\gamma>0$.

\begin{rmk}
Recall that $Y$ represents the logarithm of the (normalized) capital of a particle. Hence, the integral in \eqref{eq.game.J.def} can be interpreted as the cumulative (discounted) logarithmic utility from the dividend flow, since the dividends are typically proportional to the capital. More importantly, the dynamics of $Y$ and the objective $J$ are such that any particle will aim to maximize the drift of $Y$. The same would be true for almost any other (reasonable) objective function: e.g., we could replace $Y^x_t$ inside the integral in \eqref{eq.game.J.def} by any non-decreasing function of $Y^x_t$, covering, in particular, the case when each particle aims to maximize its expected lifetime. We view this (almost) independence of a specific choice of objective as a strength of the proposed model.
\end{rmk}

In order to define the notion of equilibrium in this model, we introduce $(\bC^b,\bC^l,\bnu)$, representing the collective controls of all particles. Following the standard methodology of continuum player games, we shall view the collective controls as a mapping from the set of particles' characteristics, $(x,y)$, to the space of individual controls. However, it turns out that, in the present case, the optimal controls of each particle depend only on its type, $x$, hence, to simplify the notation, we define: $\bC^b,\bC^l\in (L^{0}_+(\RR_+))^{|\mX|}$, and $\bnu\in \bar{\mV}$, where we introduce
$$
\bar{\mV}_0 := \left\{ \bnu: \mX\rightarrow \RR^{|\mX|}_+:\, \bnu(x,\cdot)\in \mV_0,\text{ for all $x\in\mX$}\right\},
\quad
\bar{\mV} := L^{0}_+(\RR_+\rightarrow \bar{\mV}_0).
$$
The functions $(r,\theta,\bC^b,\bC^l,\bnu)$ are determined endogenously in equilibrium.


\begin{definition}\label{def:game.optimalCont}
For a given pair $(r,\theta)\in (L^{0}_+(\RR_+))^{|\mX|}\times\Theta$ and (initial) characteristics $(x,y)\in\mX\times(0,\infty)$, the control $(C^b,C^l,\nu)$ is optimal if it is admissible, and
$$
J(x,y;r,\theta;C^b,C^l,\nu) \geq J(x,y;r,\theta;\tilde{C}^b,\tilde{C}^l,\tilde{\nu}), 
$$
for any other admissible control $(\tilde{C}^b,\tilde{C}^l,\tilde{\nu})$.
\end{definition}

\begin{definition}\label{def:game.Equil}
$(r,\theta,\bC^b,\bC^l,\bnu)\in (L^{0}_+(\RR_+))^{|\mX|}\times\Theta\times \left((L^{0}_+(\RR_+))^{|\mX|}\right)^2\times \bar{\mV}$ form an equilibrium if 
\begin{enumerate}
\item for every $(x,y)\in\mX\times(0,\infty)$, $(\bC^b(x),\bC^l(x),\bnu(x,\cdot))$ is optimal;
\item the market clears: for every $x\in\mX$ and almost every $t\geq0$,
\begin{equation}\label{eq.game.clearingCond}
\bC^b_t(x) = \int_{\mathcal{X}} \bnu_t(x',x) \bC^l_t(x') \mu(dx');
\end{equation}
\item and the model is consistent: for every $x\in\mX$ and every $t\geq0$,
\begin{equation}\label{eq.game.consistCond}
\theta_t(x) = \PP\Big(\inf_{s\in[0,t]} Y^x_s (r,\theta;\bC^b(x),\bC^l(x),\bnu(x,\cdot)) >0\Big).
\end{equation}
\end{enumerate}
\end{definition}

In the remainder of this section, we show how to construct such an equilibrium.

\subsection{Construction of equilibrium: general approach}
\label{subse:game.approach}

We will construct an equilibrium satisfying the following ansatz. 

\begin{ansatz}\label{ans:game.ans1}
There exists a default rate $\lambda\in (L^0_-(\RR_+))^{|\mX|}$ such that, for every $x\in\mX$,
$$
\theta_t(x) = 1 + \int_0^t \lambda_s(x)ds >0,\quad t\geq0.
$$
We let
$$
\blambda_t(x):=\frac{\mathrm{d}}{\mathrm{d}t}\log \theta_t(x)=\frac{\lambda_t(x)}{\theta_t(x)},\quad t\ge0.
$$
\end{ansatz}

The above ansatz implies that no cascades are possible in the equilibrium we construct. The economic explanation of this is as follows. If the distribution of particles' health levels, within a given type, approaches a critical level, at which a cascade would be triggered by this type, the lenders will simply stop lending to the particles of this type (i.e. the total weight of the in-coming edges for this node will vanish), thus, eliminating the possibility of contagion. A side effect of such action, of course, is that the particles in distress are more likely to die (as their drifts are reduced due to the lack of borrowing opportunities)\footnote{It is natural to baptize our construction a ``controlled-fire equilibrium", by analogy with the common practice in forest management. We thank Ronnie Sircar for suggesting this name.}.
In the equilibrium constructed herein, the above ansatz holds, and the particles act according to this principle.

\begin{rmk}
The above conclusions, and the entire model proposed herein, rely on the following assumptions: (i) the particles are well informed about the default intensities at each type, and (ii) they have a perfect control over their connections. In reality, of course, none of these assumptions is satisfied fully, and one expects to observe systems whose behavior lies somewhere in between the perfectly controlled and completely uncontrolled systems.
\end{rmk}

Under the above ansatz, the admissible dynamics of $Y^x$, prior to default, can be rewritten as follows:
\begin{equation}\label{eq.game.controlledYdyn}
\mathrm{d}Y^x_t = b(r_t(x),\blambda_t,C^b_t,C^l_t,\nu_t;x)\,\mathrm{d}t + \sigma\,\mathrm{d}W_t,
\end{equation}
where
\begin{equation*}\label{eq.game.GemEquil.b.def}
\begin{split}
b(r_t(x),\blambda_t,C^b_t,C^l_t,\nu_t;x) := \alpha'(x)
+ (C^b_t-C^l_t) \alpha(x)
-  C^b_t r_t(x) \qquad \\
+ C^l_{t} \int_{\mX} \nu_{t}(x') \left(r_t(x') + \blambda_t(x')\right) \mu(\mathrm{d}x').
\end{split}
\end{equation*}
Our strategy for constructing an equilibrium is based on the following observation. Note that the dynamics of $\{Y^x\}_{x\in\mX}$ are coupled across $x\in\mX$ via the interest rates $\{r(x)\}_{x\in\mX}$, the collective controls $\{\bC^b(x),\bC^l(x),\bnu(x,\cdot)\}_{x\in\mX}$, and the default rates $\{\blambda(x)\}_{x\in\mX}$. The former two affect the dynamics of individual particles, as well as the optimality, market clearance, and consistency conditions, at time $t$ only through their values at that time, i.e. the interaction is completely local. In contrast, the interaction through $\blambda$ is not completely local: the past values of $\blambda$ affect future dynamics of individual particles, and, vice versa, to find the value of $\blambda_t(x)$, one has to know the dynamics of the system up to time $t$. As a result, the equilibrium conditions for $(r,\bC^b,\bC^l,\bnu)$ can be resolved for each time $t$ separately, leading to a static fixed-point problem, which resembles the network flow problem. However, the equilibrium conditions for $\blambda$ lead to a dynamic fixed-point problem. In the subsequent construction, we separate the two problems.

To better explain our approach, let us, first, describe the optimal strategies of the particles.

\begin{lemma}\label{le:game.genEquil.opt}
Consider any $x\in\mX$, $(r,\theta)\in (L^{0}_+(\RR_+))^{|\mX|}\times\Theta$ and $(C^b,C^l,\nu)\in (L^{0}_+(\RR_+))^2\times \mV$.
Assume that $\theta$ satisfies Ansatz \ref{ans:game.ans1}, with the associated $\blambda$. Assume also that, for almost every $t\geq0$, we have:
\begin{eqnarray}
&& \nu_t \bone_{\{x'\,:\,O(x,x')=0\}}\equiv0, \label{eq.game.lemma1.admisNu} \\
&& \int_{\mX}\nu_t(x')\left(r_t(x')+\blambda_t(x')\right) \mu(\mathrm{d}x') = R_t(x), \label{eq.game.lemma1.optNu} \\
&& r_t(x) = \alpha(x)\vee R_t(x), \label{eq.game.lemma1.optr} \\
&& C^b_t \geq0 ,\quad C^l_t = (\bC + C^b_t)\bone_{\{R_t(x)>\alpha(x)\}}, \label{eq.game.lemma1.optC}
\end{eqnarray}
with 
\begin{equation}\label{eq.game.lemma1.RL}
R_t(x):= \Big(\max_{x'\,:\,O(x,x')=1}\left(r_t(x')+\blambda_t(x')\right)\!\Big)^+.
\end{equation}
If, in addition, $b(r_t(x),\blambda_t,C^b_t,C^l_t,\nu_t;x)$ is locally integrable in $t\in\RR_+$,
then $(C^b,C^l,\nu)$ is optimal for any particle with (initial) characteristics $(x,y)$, for any $y>0$. 
\end{lemma}

\begin{rmk}\label{rem:game.DomLenders.lemma1}
The conditions (\ref{eq.game.lemma1.optNu}) and (\ref{eq.game.lemma1.RL}) can be interpreted as follows.
The lenders lend everything to the neighbors that offer the best profit, i.e. have the largest $r_t+\blambda_t$ (provided it exceeds the profits from external investment), and borrowing generates zero profits. This, in particular, implies that $\nu_t$ must be supported on the set on which this maximum is attained, and it is consistent with the implicit assumption of ``dominant lenders" (or ``competitive borrowers"), as discussed in Remark \ref{rem:game.CompetitiveBorrowers}. Of course, this is an extreme case -- in practice, borrowing should result in some profit -- but, as explained below, this assumption allows for a tractable equilibrium, and more realistic equilibria can be constructed as perturbations of the proposed one. 
\end{rmk}

\noindent\textbf{Proof.} It is easy to see that $(C^b,C^l,\nu)$ is admissible. Consider any other admissible control $(\tilde{C}^b,\tilde{C}^l,\tilde{\nu})$ and denote by $b$ and $\tilde{b}$ the associated drifts. Notice that (\ref{eq.game.lemma1.optNu}) and (\ref{eq.game.lemma1.RL}) imply:
$$
\int_{\mX} \tilde{\nu}_{t}(x') \left(r_t(x') + \blambda_t(x')\right) \mu(\mathrm{d}x') \leq R_t(x)
= \int_{\mX}\nu_t(x')\left(r_t(x')+\blambda_t(x')\right) \mu(\mathrm{d}x').
$$
In addition, (\ref{eq.game.lemma1.optr}) and (\ref{eq.game.lemma1.optC}) yield
\begin{equation*}
(C^b_t,C^l_t) \in \text{argmax}_{(y,z)\,:\,y\geq 0,\,0\leq z\leq \bC + y}\left( y\left(\alpha(x)-r_t(x)\right) + z \left(R_t(x) - \alpha(x)\right)\right).
\end{equation*}
Recalling (\ref{eq.game.GemEquil.b.def}) we conclude that $b\geq \tilde{b}$ and, hence, $J(x,y;r,\theta;C^b,C^l,\nu) \geq J(x,y;r,\theta;\tilde{C}^b,\tilde{C}^l,\tilde{\nu})$. \ep

\medskip

Extending the above result in an obvious way, we obtain a sufficient condition for equilibrium.

\begin{cor}\label{cor:game.EquilSuffCond.Dyn}
Consider any $(r,\theta,\bC^b,\bC^l,\bnu)\in (L^{0}_+(\RR_+))^{|\mX|}\times\Theta\times \left((L^{0}_+(\RR_+))^{|\mX|}\right)^2\times \bar{\mV}$.
Assume that, for every $x\in\mX$, the quadruple $(r,\theta,\bC^l,\bnu)$ and the control $(C^b:=\bC^b(x),C^l:=\bC^l(x),\nu:=\bnu(x,\cdot))$ satisfy the assumptions of Lemma \ref{le:game.genEquil.opt}. If, in addition, (\ref{eq.game.clearingCond}) holds for almost every $t$ and every $x$, and (\ref{eq.game.consistCond}) holds for every $t$ and $x$, then $(r,\theta,\bC^b,\bC^l,\bnu)$ form an equilibrium.
\end{cor}

To explain our solution approach, we point out that, for given $t$ and $\blambda_t(\cdot)$, the conditions (\ref{eq.game.clearingCond}) and (\ref{eq.game.lemma1.admisNu})--(\ref{eq.game.lemma1.RL}) can be viewed as a system of equations for $\{r_t(x),\bC^b_t(x),\bC^l_t(x),\bnu_t(x)\}_{x\in\mX}$. Assume that we can solve the static problem, i.e. we can find a function that maps any $\blambda_t(\cdot)$ to a solution of this system. Then, to construct an equilibrium, we only need to solve the dynamic problem, i.e. find a fixed-point of the mapping from $\blambda$ to the time-derivative of the logarithm of the right-hand side of (\ref{eq.game.consistCond}), computed via the solution to the static problem $(r(\blambda),\bC^b(\blambda),\bC^l(\blambda),\bnu(\blambda))$. The solutions to the two problems are described in the remainder of this section.

\subsection{Solving the dynamic problem}
\label{subse:game.dynamic}

In this subsection, we assume that the static problem has been solved (in the next subsection, we construct such a solution).
\begin{ass}\label{ass:game.1}
Assume that there exists a measurable mapping
$$
\mS\,:\, \RR^{|\mX|}_-\ni\blambda_t \mapsto \big(\hat{r}(\blambda_t),\hC^b(\blambda_t),\hC^l(\blambda_t),\hnu(\blambda_t)\big)\in \big(\RR^{|\mX|}_+\big)^3\times \bar{\mV}_0
$$
such that,
\begin{enumerate}
\item for every $\blambda_t\in \RR^{|\mX|}_-$ and every $x\in\mX$, 
\begin{equation*}
\begin{split}
r_t:=\hat{r}(\blambda_t),\,\,\bC^b_t:=\hC^b(\blambda_t),\,\,\bC^l_t:=\hC^l(\blambda_t),\,\,\bnu_t:=\hnu(\blambda_t),
C^b_t:=\bC^b_t(x), \,\,C^l_t:=\bC^l_t(x),\\ 
\nu_t:=\bnu_t(x,\cdot)
\end{split}
\end{equation*}
satisfy (\ref{eq.game.clearingCond}) and (\ref{eq.game.lemma1.admisNu})--(\ref{eq.game.lemma1.RL});
\item the mapping 
$$
\RR^{|\mX|}_-\ni\blambda_t\mapsto \Big(x\mapsto \hb(\blambda_t;x) := b\big(\hr(\blambda_t;x),\blambda_t,\hC^b(\blambda_t;x),\hC^l(\blambda_t;x),\hnu(\blambda_t;x,\cdot);x\big)\Big)
\in \RR^{|\mX|},
$$
with $b$ given by (\ref{eq.game.GemEquil.b.def}), is Lipschitz and bounded absolutely by a constant.
\end{enumerate}
\end{ass}

In the next subsection, we construct $\mS$ satisfying the above assumption. It is worth mentioning that the construction of such $\mS$ is the main mathematical challenge resolved in the present section.
In this subsection, we assume that the desired $\mS$ is given.

Notice that, for any measurable $[0,T]\ni t\mapsto \blambda_t\in \RR^{|\mX|}_-$, the mappings 
$$
(t,x,x')\mapsto \hat{r}(\blambda_t;x),\,\hC^b(\blambda_t;x),\,\hC^l(\blambda_t;x),\,\hnu(\blambda_t;x,x')
$$
are measurable, and $(r,\bC^b,\bC^l,\bnu)$, defined in part 1 of Assumption \ref{ass:game.1}, satisfies all assumptions of Corollary \ref{cor:game.EquilSuffCond.Dyn}, except, possibly, (\ref{eq.game.consistCond}).

The main goal of this subsection is to construct $\blambda$ such that the associated $(r,\theta,\bC^b,\bC^l,\bnu)$ satisfies (\ref{eq.game.consistCond}). To this end, we fix an arbitrary $T\in(0,\infty)$ and analyze the fixed-point problem of the mapping $\mD:\,\mathbb{L}^{2}\big([0,T]\rightarrow\RR^{|\mX|}_-\big)\to \mathbb{L}^0\big([0,T]\rightarrow \RR^{|\mX|}_-\big)$ given by
\begin{equation}\label{eq.mD.def}
\mD:\,\blambda\mapsto 
\bigg((t,x)\mapsto\frac{\mathrm{d}}{\mathrm{d}t} \log \PP\Big(\inf_{s\in[0,t]} Y^x_s(\blambda) >0\Big)\bigg),
\end{equation}
where
\begin{equation}\label{eq.YofBlambda.def}
Y^x_t(\blambda) := Y^x_0 + \int_0^t \hb(\blambda_s;x)\,\mathrm{d}s + \alpha'(x) t + \sigma W_t.
\end{equation}
Note that the mapping $\mD$ is well-defined due to Assumption \ref{ass:game.1}. Indeed, for any measurable $[0,T]\ni t\mapsto \blambda_t\in \RR^{|\mX|}_-$, the mapping $(t,x)\mapsto \hb(\blambda_t;x)$ is measurable. The latter and the boundedness of $\hb$ imply that $Y^x(\blambda)$ is well-defined.
In addition, for any two modifications $\blambda^1$ and $\blambda^2$, for every $x$, we have: $\hb(\blambda^1_t;x)=\hb(\blambda^2_t;x)$ for almost every $t$, and, in turn, $Y^x(\blambda^1)_t= Y^x_t(\blambda^2)$ for all $t$.

Our analysis of the fixed-point problem for $\mathcal{D}$ is based on the associated system of PDEs for $p=\{p(\cdot,\cdot;x)\}_{x\in\mX}$, where $p(t,\cdot;x)$ is the density of the restriction of the law of $Y^x_t$ to $(0,\infty)$. In accordance with the existing results (cf. \cite{NadShk}, \cite{Delarue1}), we expect $p$ to satisfy the following system of equations:
\begin{equation}\label{eq.game.PDE.1} 
p_t(t,y;x) = -\hb_t(\blambda_t;x)\,p_y(t,y;x) + \frac{1}{2}\sigma^2 p_{yy}(t,y;x),\quad p(0,y;x) = p_0(y;x),\quad p(t,0;x)=0, 
\end{equation}
\begin{equation}
\blambda_t(x) = -\frac{\sigma^2}{2} \frac{p_y(t,0;x)}{\int_0^{\infty} p(t,y;x)\,\mathrm{d}y}. \label{eq.game.PDE.2}
\end{equation}
The following proposition proves that there exists a unique solution to the above system, in turn, proving the existence and uniqueness of the fixed-point of $\mD$. It is based on the results of \cite{NadShk}.
In this proposition, and in its proof, we denote by $W^n_2(\RR_+)$ and $W^{n,m}_2([0,T]\times\RR_+)$ the standard Sobolev spaces.

\begin{proposition}\label{prop:game.prop1}
Consider a mapping $\mS$ satisfying Assumption \ref{ass:game.1} and assume that, for every $x$, $Y^x_0\sim p_0(y;x)\,\mathrm{d}y$ with $p_0(\cdot;x)\in W^1_2(\RR_+)$ and $p_0(0;x)=0$. Then, there exists a unique fixed-point of the associated mapping $\mD$.
\end{proposition}

\noindent\textbf{Proof.}
Let us fix an arbitrary $\blambda\in\mathbb{L}^2\big([0,T]\rightarrow\RR^{|\mX|}_-\big)$.
Then, repeating steps 1 and 2 in the proof of \cite[Lemma 3.1]{NadShk}, we prove the following statement. For every $x\in\mX$, we have that, for all $t\in[0,T]$, the distribution of 
$$
Y^x_t(\blambda) \bone_{\{\inf_{s\in[0,t]} Y^x_s(\blambda) >0\}},
$$
restricted to $(0,\infty)$, has a density $p(t,\cdot;x)$ (recall that $Y^x(\blambda)$ is defined in (\ref{eq.YofBlambda.def})). Moreover, $p(\cdot,\cdot;x)$ is the unique solution of (\ref{eq.game.PDE.1}) in $W^{1,2}_2([0,T]\times\RR_+)$.
Lemma 3.2 of \cite{NadShk} shows that, for every $x$ and almost every $t$, $\mD[\blambda](t,x)$ is given by the right hand side of (\ref{eq.game.PDE.2}). In addition, the last paragraph of step 1 in the proof of \cite[Proposition 4.1]{NadShk} shows that, for every $x$, $\mD[\blambda](\cdot,x)\in\mathbb{L}^2([0,T])$, hence, the range of $\mD$ is contained in $\mathbb{L}^2([0,T]\rightarrow \RR^{|\mX|}_-)$.

Next, step 3 in the proof of \cite[Proposition 4.1]{NadShk} shows that there exists an increasing function $C_1:[0,T]\rightarrow\RR_+$, converging to $0$ at $0$, such that, for any $\blambda^1,\blambda^2\in \mathbb{L}^2([0,T]\rightarrow\RR^{|\mX|}_-)$, we have, for every $x$:
$$
\int_0^{s} \left(\mD[\blambda^1](t,x) - \mD[\blambda^2](t,x)\right)^2\,\mathrm{d}t
\leq C^2_1(s) \int_0^{s} \left(\blambda^1(t,x) - \blambda^2(t,x)\right)^2 \,\mathrm{d}t,\quad s\in[0,T],
$$
which yields the contraction property of $\mD$ for sufficiently small $T>0$.
Then, the Banach fixed-point theorem yields the existence and uniqueness of the fixed-point of $\mD$, for sufficiently small $T>0$. Step 4 in the proof of \cite[Proposition 4.1]{NadShk} shows how to extend this result to arbitrary $T\in(0,\infty)$. \ep

\medskip

Note that the domain of $\mD$ can be extended in a straightforward way to the space of locally square integrable functions $\RR_+\rightarrow\RR^{|\mX|}_-$.
Proposition \ref{prop:game.prop1} then implies that there exists a unique fixed-point of $\mD$ in this space.
Whenever we refer to a fixed-point $\blambda$ of $\mD$, we understand it as an element of this space.

For any $\blambda\in L^0(\RR_+\rightarrow\RR^{|\mX|}_-)$, we let
$$
\hat{\theta}_t(\blambda;x) := \exp\left( \int_0^t \blambda_s(x) ds\right).
$$
Notice that, if $\blambda$ is (a modification of) the fixed-point of $\mD$, then $\theta: (t,x)\mapsto \theta_t(x):= \hat{\theta}_t(\blambda;x)$ belongs to $\Theta$, and $(r,\theta,\bC^b,\bC^l,\bnu)$, with $(r,\bC^b,\bC^l,\bnu)$ defined in part 2 of Assumption \ref{ass:game.1}, satisfies all assumptions of Corollary \ref{cor:game.EquilSuffCond.Dyn}, including (\ref{eq.game.consistCond}).
Thus, we have proved the following result.

\begin{cor}\label{cor:game.EquilSuffCond.Dyn.2}
Consider a mapping $\mS$, satisfying Assumption \ref{ass:game.1}, with the associated $(\hat{r},\hC^b,\hC^l,\hnu)$. Assume that, for every $x\in\mX$, $Y^x_0\sim p_0(y;x)\,\mathrm{d}y$ with $p_0(\cdot;x)\in W^1_2(\RR_+)$ and $p_0(0;x)=0$. Let $\blambda$ be the fixed-point of the associated mapping $\mD$, and let
\begin{equation*}
\begin{split}
r_t(x):=\hat{r}(\blambda_t;x),\,\,\theta_t(x):= \hat{\theta}_t(\blambda;x),
\,\,\bC^b_t(x):=\hC^b(\blambda_t;x),\,\,\bC^l_t(x):=\hC^l(\blambda_t;x), \\
 \bnu_t(x,x'):=\hnu(\blambda_t;x,x').
\end{split}
\end{equation*}
Then, $(r,\theta,\bC^b,\bC^l,\bnu)$ form an equilibrium.
\end{cor}

\subsection{Solving the static problem}
\label{subse:game.static}


Herein, we present a solution to the static fixed-point problem, i.e. we construct a mapping $\mS$ satisfying Assumption \ref{ass:game.1}. Consistent with the initial formulation of the problem and with the equations (\ref{eq.game.lemma1.optNu})--(\ref{eq.game.lemma1.optr}) (cf. Remarks \ref{rem:game.CompetitiveBorrowers} and \ref{rem:game.DomLenders.lemma1}), the specific choice of the construction presented can be interpreted as the assumption of ``dominant lenders" (or ``competitive borrowers").

\begin{rmk}
Even under the implicit assumption of dominant lenders, the uniqueness of the equilibrium is not clear. Nevertheless, certain characteristics of equilibria (such as the interest rate) may be unique, as discussed in Remark \ref{rem:game.uniqueness}.
In any case, the (constructive) proof of Proposition \ref{prop:game.prop2} shows that there exists a natural choice of such an equilibrium.
\end{rmk}

As the derivations in this subsection are performed for a fixed time $t$, for convenience, we drop this subscript. Below we formulate a sufficient condition for $\mS:\blambda\mapsto (\hr,\hC^b,\hC^l,\hnu)$ to satisfy Assumption \ref{ass:game.1} (i.e. to be a solution of the static problem).
Assume that we are given $(r,\blambda,\bC^b,\bC^l,\bnu)\in (L^{0}_+(\RR_+))^{|\mX|}\times (L^{0}_-(\RR_+))^{|\mX|}\times \left(L^{0}_+(\RR_+))^{|\mX|}\right)^2\times \bar{\mV}$ and consider the following conditions.

\begin{enumerate}
\item $r$ satisfies
\begin{equation}\label{eq.r.equil.def}
r(x) = \alpha(x)\vee \max_{y\,:\,O(x,y)=1}(r(y)+\blambda(y)),\quad x\in\mX.
\end{equation}
\item For every $x\in\mX$,
\begin{equation}\label{eq.nu.equil.def}
r(x)>\alpha(x)\,\,\Rightarrow\,\,
\bnu(x,x')=\bone_{A(x)}(x')/\mu(A(x)),\quad x'\in\mX,
\end{equation}
where $A(x)$ is a non-empty subset of
\begin{equation}\label{eq.game.tildeA.def}
\tilde{A}(x):= \{z\in\mX:\,O(x,z)=1,\, r(z)+\blambda(z) = \max_{y\,:\,O(x,y)=1}(r(y)+\blambda(y))\},
\end{equation}
if $r(x)\leq \alpha(x)$, then $\bnu(x,\cdot)\equiv0$.

\item For all $x\in\mX$, $\bC^b$and $\bC^l$ satisfy
\begin{equation}
\label{eq.Cbl.equil.def}\bC^b(x) = \int_{\mathcal{X}} \bnu(x',x) \,\bC^l(x') \mu(dx'),\quad 
\bC^l(x) = (\bC + \bC^b(x))\,\bone_{\{R(x)>\alpha(x)\}},
\end{equation}
where we recall (cf. (\ref{eq.game.lemma1.RL}))
\begin{equation}\label{eq.game.lemma1.R.2}
R(x)= \Big(\max_{x'\,:\,O(x,x')=1}(r(x')+\blambda(x'))\Big)^+.
\end{equation}
\end{enumerate}
To explain the rationale behind equations (\ref{eq.r.equil.def})--(\ref{eq.Cbl.equil.def}), we recall the principle discussed in Remarks \ref{rem:game.CompetitiveBorrowers} and \ref{rem:game.DomLenders.lemma1} (compare (\ref{eq.r.equil.def})--(\ref{eq.game.lemma1.R.2}) to (\ref{eq.game.lemma1.admisNu})--(\ref{eq.game.lemma1.RL})).
Namely, in the equilibrium we construct, the particles lend everything to the neighbors that offer the best profit, i.e. have the largest $r_t+\blambda_t$, provided it exceeds their profits from external investment. Hence, $\nu_t$ must be supported on the set on which this maximum is attained. Finally, we use the positive part operator in the expression for $R$, as the latter is only relevant when it exceeds $\alpha$, which is nonnegative.

Notice that, once $r$ and $\blambda$ are fixed, the above equations determine the network flow $(\bC^b,\bC^l,\bnu)$ that maximizes the expected profit of each member of the network, provided the other members do not deviate from these strategies. This is precisely the analogue of the network flow problem in the presence of multiple agents (i.e. when the network optimization is decentralized).

Notice that, if (\ref{eq.r.equil.def})--(\ref{eq.Cbl.equil.def}) are satisfied, we have:
$$
r(x) = \alpha(x)\vee R(x),
$$
and the associated drift, defined in (\ref{eq.game.GemEquil.b.def}), becomes
\begin{equation}\label{eq.game.b.Static}
\begin{split}
b(x) &= b(r(x),\blambda,\bC^b(x),\bC^l(x),\bnu(x,\cdot);x) \\
&= \alpha'(x) 
+ \bC^b(x) (\alpha(x)-r(x)) + \bC^l(x) \left( \int_{\mathcal{X}} \bnu(x,x')(r(x')+\blambda(x')) \mu(\mathrm{d}x') -\alpha(x)\right) \\
&= \alpha'(x) + \bC \left(R(x)-\alpha(x)\right)^+.
\end{split}
\end{equation}

It is easy to notice that conditions (\ref{eq.r.equil.def})--(\ref{eq.game.lemma1.R.2}) are a special case of (\ref{eq.game.clearingCond}), (\ref{eq.game.lemma1.admisNu})--(\ref{eq.game.lemma1.RL}). This observation is formalized in the following lemma.

\begin{lemma}\label{le:game.staticSuffCond}
Let $(r,\blambda,\bC^b,\bC^l,\bnu)\in (L^{0}_+(\RR_+))^{|\mX|}\times (L^{0}_-(\RR_+))^{|\mX|}\times \left((L^{0}_+(\RR_+))^{|\mX|}\right)^2\times \bar{\mV}$ satisfy (\ref{eq.r.equil.def})--(\ref{eq.Cbl.equil.def}).
Then, for every $x$,
$$
\left(r_t:=r,\,\,\bC^b_t:=\bC^b,\,\,\bC^l_t:=\bC^l_t,\,\,\bnu,\,\,C^b_t:=\bC^b(x), \,\,C^l_t:=\bC^l(x),\,\,\nu_t:=\bnu(x,\cdot)\right)
$$
satisfy (\ref{eq.game.clearingCond}) and (\ref{eq.game.lemma1.admisNu})--(\ref{eq.game.lemma1.RL}).
\end{lemma}

\noindent\textbf{Proof.} The lemma follows by direct verification. \ep

\medskip

Finally, we present the main mathematical result of this section, which proves the existence of a solution to (\ref{eq.r.equil.def})--(\ref{eq.Cbl.equil.def}) and its regularity with respect to $\blambda$.

\begin{proposition}\label{prop:game.prop2}
There exists a mapping 
$$
\mS:\, \RR_-^{|\mX|}\ni\blambda\mapsto \big(\hat{r}(\blambda), \hat{C}^b(\blambda), \hat{C}^l(\blambda), \hat{\nu}(\blambda)\big)\in
\big(\RR_+^{|\mX|}\big)^3\times \bar{\mathcal{V}}_0
$$
such that the associated 
$$
\big(r:=\hat{r}(\blambda),\,\blambda,\, \bC^b:=\hat{C}^b(\blambda),\, \bC^l:=\hat{C}^l(\blambda),\, \bnu:=\hat{\nu}(\blambda)\big)
$$
satisfy (\ref{eq.r.equil.def}), (\ref{eq.nu.equil.def}), (\ref{eq.Cbl.equil.def}), for every $\blambda$. Moreover, for every $x\in\mX$, the associated $R(x)$, given by (\ref{eq.game.lemma1.R.2}), is Lipschitz and bounded in $\blambda$.
\end{proposition}


\noindent\textbf{Proof.} For any given $\blambda$, we need to construct $(r,\bC^b,\bC^l,\bnu)$ satisfying (\ref{eq.r.equil.def}), (\ref{eq.nu.equil.def}), (\ref{eq.Cbl.equil.def}). Without loss of generality we assume $\mX=\{1,\ldots,n\}$.

\medskip

\noindent\textbf{Step 1: solving (\ref{eq.r.equil.def}).} Notice that the equation (\ref{eq.r.equil.def}) is nothing else but a system of linear equations in the max-plus algebra. Recall that the max-plus algebra is a semiring on $\RR\cup\{-\infty\}$, in which the addition is replaced by the maximum operator, denoted by $\oplus$, and the multiplication is replaced by the addition, $\otimes$:
$$
x\oplus y = \max(x,y),\quad x\otimes y = x+y,\quad x,y\in \RR\cup\{-\infty\}.
$$
Note that the matrix addition and multiplication in the max-plus algebra can be defined in the usual way, via the above operations, and they satisfy the usual associativity and distributivity rules. Then, (\ref{eq.r.equil.def}) reads as
$$
r = A\otimes r \oplus \alpha,\quad r,\alpha\in (\RR\cup\{-\infty\})^n,\,\,A\in (\RR\cup\{-\infty\})^{n\times n},
$$
where $\alpha_i = \alpha(i)$, for $i=1,\ldots,n$, and $A$ is a square matrix, with
$$
A_{ij}= \blambda(j)\,\bone_{\{O(i,j)=1\}} - \infty\cdot \bone_{\{O(i,j)=0\}},\quad i,j=1,\ldots,n
$$
and the convention $\infty\cdot 0 = 0$.
Since the entries of $A$ are nonpositive, \cite[Theorems 3.20, 4.75, and Remark 4.80]{Baccelli} imply that
\begin{equation}\label{eq.game.minsol}
r: = A^*\otimes \alpha, 
\end{equation}
where the Kleene star matrix
$$
A^*:=I \oplus A \oplus A^2 \oplus \cdots\oplus A^{n-1},\quad I_{ij}= -\infty\cdot \bone_{\{i\neq j\}},
$$
is the smallest solution to (\ref{eq.r.equil.def}), in the sense that it solves (\ref{eq.r.equil.def}), and $r_i\leq \tilde{r}_i$, for all $i$, for any solution $\tilde{r}$.
Note that the associated vector $R$, with $R_i:=R(i)$, is given by
$$
R = A\otimes r \oplus e = (A \oplus A^2 \oplus \cdots\oplus A^{n})\otimes \alpha \oplus e,
$$
where $e$ is the vector of zeros. It is clear that $R_i \leq \max_j \alpha_j $, for all $i$.
In addition, it is easy to see that, if $v$ is a Lipschitz function of $\blambda$, with values in $\RR^n$, then $A\otimes v$ is also Lipschitz in $\blambda$. The operator $\oplus$ preserves the Lipschitz property as well. Thus, by induction, $R$ is Lipschitz in $\blambda$.

\medskip

\noindent\textbf{Step 2: eliminating cycles.}
Let us show that there exists an optimal lending allocation, denoted by $\tilde{\nu}^m$, that satisfies (\ref{eq.nu.equil.def}) and produces no lending cycles.
To this end, we define $\tilde{\nu}\in \bar{\mV}$:
\begin{equation*}
\tilde{\nu}(i,j)=
\left\{
\begin{array}{ll}
{\bone_{\tilde{A}(i)}(j)/\mu(\tilde{A}(i)),} & {r(i)>\alpha(i),}\\
{0,} & {r(i)=\alpha(i),}
\end{array}
\right.
\end{equation*}
where $\tilde{A}(i)$ is given by (\ref{eq.game.tildeA.def}).
Note that $\tilde{\nu}$ satisfies (\ref{eq.nu.equil.def}).
Next, consider the new adjacency matrix $\tilde{O}$ induced by $\tilde{\nu}$,
$$
\tilde{O}(i,j) = \bone_{\{\tilde{\nu}(i,j)>0\}},
$$
and the associated directed graph. In the new graph $(\mX,\tilde{O})$, let us pick a node $x_0\in\mX$ and denote by $C$ the union of all cycles containing this node. We assume that there exists an $x_0\in\mX$ such that $C$ is non-empty (otherwise, $(\mX,\tilde{O})$ contains no cycles, and we can skip this step). We let $C=(x_0,x_1,\ldots,x_{k-1})$, with $k\geq1$. The set $C$ has the following properties: (i) $C$ is connected, (ii) any two $x,y\in C$ are contained in at least one cycle, and (iii) all cycles containing at least one $x\in C$ are included in $C$.

Note that, for every cycle $(x_{i_0},\ldots,x_{i_{l-1}})\subset C$, we have 
$$
R(x_{i_j})=r(x_{i_j}) = \text{max}_{y\,:\,O(x_{i_j},y)=1}(r(y)+\blambda(y)) = r(x_{i_{(j+1)\text{ mod }l}}) + \blambda(x_{i_{(j+1)\text{ mod }l}}).
$$
Since $\blambda\leq 0$ and any two points in $C$ are contained in a cycle in $C$, the above implies 
$$
\blambda(x_0)=\cdots=\blambda(x_{k-1})=0,
\quad r(x_0)=\cdots=r(x_{k-1})=:r^c>\text{max}_{x_i\in C} \,\alpha(x_i) \geq0.
$$
Next, for every $x_i\in C$, we define the ``alternative profit"
$$
R^a(x_i):= \text{max}_{y\,:\,O(x_i,y)=1,\,y\notin C}(r(y)+\blambda(y)) \leq r^c.
$$

Our next goal is to show that there exists an $x_i\in C$ such that $R^a(x_i) = r^c$.
Assume the opposite: $R^a(x_i)< r^c$, for every $i=0,1,\ldots,k-1$. Consider a new graph $(\hat{\mX},\hat{O})$, which only differs form $(\mX,\tilde{O})$ in that the set $C$ is ``collapsed'' into a single state, denoted by $\bar{x}$. Setting $\alpha(\bar{x})=r^c$, $\blambda(\bar{x})=0$, and $r(\bar{x})=r^c$, we note that $r$ solves (\ref{eq.r.equil.def}) for the new network $(\hat{\mX},\hat{O},\alpha)$.
On the other hand, choosing a constant
$$
\hat{r}^c\in \left( \text{max}_{i=0,\ldots,k-1}\,R^a(x_i)\vee\alpha(x_i),\, r^c\right)
$$
and setting $\hat{\alpha}(\bar{x}):=\hat{r}^c$ (with the other values of $\hat{\alpha}$ being the same as those of $\alpha$), we can find the minimal solution $\hat{r}$ to (\ref{eq.r.equil.def}) for the network $(\hat{\mX},\hat{O},\hat{\alpha})$. It is clear from (\ref{eq.r.equil.def}) that $\hat{r}(\bar{x})\geq\hat{\alpha}(\bar{x})=\hat{r}^c$. In addition, from the minimality of $\hat{r}$ and the monotonicity of the right hand side of (\ref{eq.game.minsol}) with respect to $\alpha$, we conclude that $\hat{r}(x)\leq r(x)$, for every $x\in\hat{O}$. This observation, along with our assumption and the choice of $\hat{r}^c$, imply that
$$
\text{max}_{y\,:\,\hat{O}(\bar{x},y)=1}(\hat{r}(y)+\blambda(y)) < \hat{r}^c=\hat{\alpha}(\bar{x}).
$$
This, in turn, yields that $\hat{r}(\bar{x})=\hat{\alpha}(\bar{x})=\hat{r}^c<r^c$. Finally, we turn back to the network $(\mX,\tilde{O},\alpha)$ and, expanding $\hat{r}$ in an obvious way (i.e. setting $\hat{r}(x)=\hat{r}^c$ for all $x\in C$), we obtain an $\hat{r}$ that solves (\ref{eq.r.equil.def}) for $(\mX,\tilde{O},\alpha)$. Since $\hat{r}=\hat{r}^c<r^c=r$ on $C$, we arrive at a contradiction to the minimality of $r$.

Thus, there must exist an $x^1\in C$ such that $R^a(x^1) = r^c$. The latter equality implies that $\tilde{A}(x^1)$ contains at least one element outside of $C$.
Then, we can modify $\tilde{\nu}$ by creating
$$
\tilde{\nu}^1(x,y): =
\left\{
\begin{array}{ll}
{0,} & {x=x^1,\,y\in C,}\\
{\bone_{\tilde{A}(x)\setminus C}(y)/\mu(\tilde{A}(x)\setminus C),} & {x=x^1,\,y\notin C,}\\
{\tilde{\nu}(x,y),} & {\text{otherwise.}}
\end{array}
\right.
$$
Note that $\tilde{\nu}^1$ satisfies (\ref{eq.nu.equil.def}). Thus, we define $\tilde{O}^1$ via $\tilde{\nu}^1$, and construct the new set $C^1$ in the same way as $C$ was constructed above, but in the new graph $(\mX,\tilde{O}^1)$. Note that any cycle in $(\mX,\tilde{O}^1)$ is also a cycle in $(\mX,\tilde{O})$. In addition, the node $x^1$ is contained in at least one cycle in $(\mX,\tilde{O})$ and is not contained in any cycle in $(\mX,\tilde{O}^1)$. These observations imply that: $x^1\notin C^1$, $\tilde{\nu}^1(x,\cdot)\equiv\tilde{\nu}(x,\cdot)$ for all $x\in C^1$, and the number of cycles in $(\mX,\tilde{O}^1)$ is strictly smaller than the number of cycles in $(\mX,\tilde{O})$. Thus, the above arguments can be applied inductively, terminating in a finite number steps, when the resulting graph $(\mX,\tilde{O}^m)$ contains no cycles. The resulting $\tilde{\nu}^m$ satisfies (\ref{eq.nu.equil.def}).

\medskip

\noindent\textbf{Step 3: solving both (\ref{eq.nu.equil.def}) and (\ref{eq.Cbl.equil.def}).}
Finally, let us construct $(\bC^b,\bC^l,\bnu)$ satisfying (\ref{eq.nu.equil.def}) and (\ref{eq.Cbl.equil.def}).
As $(\mX,\tilde{O}^m)$ has no cycles, there must exist at least one node $x\in\mX$ that has no out-going edges (i.e. does not lend within the network), so that
$\tilde{\nu}^m(x,\cdot)\equiv0$. Denote the collection of all such nodes by $\mX^0$. We define
$$
\bC^{b,0}(x):=0,\,\, \bC^{l,0}:=\bC\,\bone_{\{R(x)>\alpha(x)\}}(=0),\,\, \bnu(x,\cdot):=\tilde{\nu}^m(x,\cdot),\quad x\in \mX^0.
$$
Since the nodes within $\mX^0$ are not directly connected to each other, it is easy to notice that $(\bC^{b,0},\bC^{l,0},\bnu)$ satisfy (\ref{eq.nu.equil.def}) and (\ref{eq.Cbl.equil.def}), with $\mX$ replaced by $\mX^0$.
Next, we define
$$
\mX^1:=\left\{ x\in\mX\,:\,\exists\,y\in\mX^0\text{ such that } \tilde{\nu}^m(x,y)>0 \right\},
$$
and, for every $x\in\mX^1$, we set $\bnu(x,\cdot)$ to be proportional to $\tilde{\nu}^m(x,\cdot)$ on $\mX^0$, and to be zero on $\mX\setminus\mX^0$. We define
$$
\bC^{b,1}(x)=0,\,\, \bC^{l,1}(x)=\bC,\quad x\in\mX^1
$$
and update
$$
\bC^{b,1}(x):= \sum_{y\in\mX^1} \bnu(y,x) \bC^{l,1}(y) \mu(\{y\}),\,\,
\bC^{l,1}:=(\bC + \bC^{b,1}(x))\,\bone_{\{R(x)>\alpha(x)\}},\quad x\in \mX^0.
$$
Since the nodes within $\mX^j$ are not directly connected to each other, for $j=0,1$, it is easy to see that $\mX^0\cap\mX^1=\emptyset$ and that $(\bC^{b,1},\bC^{l,1},\bnu)$ satisfy (\ref{eq.nu.equil.def}) and (\ref{eq.Cbl.equil.def}), with $\mX$ replaced by $\mX^0\cup\mX^1$.
We repeat this procedure iteratively. Namely, in the $i$-th step, we define
$$
\mX^i:=\bigg\{ x\in\mX\setminus \bigcup_{j=0}^{i-1}\mX^j\,:\,\exists\,y\in\bigcup_{j=0}^{i-1}\mX^j\text{ such that } \tilde{\nu}^m(x,y)>0 \bigg\},
$$
for every $x\in\mX^i$, we set $\bnu(x,\cdot)$ to be proportional to $\tilde{\nu}^m(x,\cdot)$ on $\bigcup_{j=0}^{i-1}\mX^j$, and to be zero on $\mX\setminus\bigcup_{j=0}^{i-1}\mX^j$, define
$$
\bC^{b,i}(x)=0,\,\, \bC^{l,i}(x)=\bC,\quad x\in\mX^i,
$$
and update, iterating over $k=i,i-1,\ldots,1$,
$$
\bC^{b,i}(x):= \sum_{y\in\mX^k} \bnu(y,x) \bC^{l,i}(y) \mu(\{y\}),\,\,
\bC^{l,i}:=(\bC + \bC^{b,i}(x))\,\bone_{\{R(x)>\alpha(x)\}},\quad x\in \mX^{k-1}.
$$
It is easy to see, by induction, that, for every $i$, the nodes within $\mX^j$ are not directly connected to each other, for $j=0,1,\ldots,i$, that $(\bigcup_{j=0}^{i-1}\mX^j)\cap\mX^i=\emptyset$, and that $(\bC^{b,i},\bC^{l,i},\bnu)$ satisfy (\ref{eq.nu.equil.def}) and (\ref{eq.Cbl.equil.def}), with $\mX$ replaced by $\bigcup_{j=0}^{i-1}\mX^j$.
The construction stops after the $p$-th step, where $p$ is the smallest integer such that there are no nodes in $\mX\setminus \bigcup_{j=0}^{p}\mX^j$ directly connected to $\bigcup_{j=0}^{p}\mX^j$, according to $\tilde{\nu}^m$.
As there are no cycles in $(\mX,\tilde{O}^m)$, we must have: $\mX=\bigcup_{j=0}^{p}\mX^j$.
Letting $\bC^b:=\bC^{b,p}$ and $\bC^l:=\bC^{l,p}$ we complete the proof. \ep

\begin{rmk}\label{rem:game.uniqueness}
Theorem 3.17 in \cite{Baccelli} shows that the solution to (\ref{eq.r.equil.def}) is unique if all entries of $\blambda$ are strictly negative, which is expected to hold at almost every time $t$.\footnote{This fact can be proven easily by comparing $Y^x$ to a Brownian motion with drift and noticing that the drift $b$, given by (\ref{eq.game.GemEquil.b.def}), is bounded absolutely in equilibrium.} Therefore, typically, one may expect to have a unique interest rate $r$ in all ``dominant lenders" (or ``competitive borrowers") equilibria, even though the sizes and allocations of credit exposures may not be unique.
\end{rmk}

\begin{rmk}
Note that it is important to find the smallest solution $r$ to the max-plus linear system (\ref{eq.r.equil.def}), in order to avoid cycles and, in turn, to satisfy the market clearing condition.
To see why, consider, for example, $\mX=\{1,2,3\}$, $\alpha\equiv0$, $\blambda\equiv0$, $O(1,2)=O(2,3)=O(3,1)=1$, with all the other entries of $O$ being zero. Note that $r(1)=r(2)=r(3)=z$, with any constant $z\geq0$, solve (\ref{eq.r.equil.def}). However, if $z>0$, it is easy to see that the market clearing condition (\ref{eq.Cbl.equil.def}) cannot be satisfied. Indeed, this situation reminds of a credit bubble, in which no particle can make positive returns on external investment, but the interest rate is strictly positive, hence, each particle aims to lend to the ones it is connected to. Since the network is a circle, the initially lent funds reach the creditor, thus encouraging him to lend more. It is clear that this problem can be easily fixed by choosing $z=0$, which is the smallest solution, and also the only economically meaningful solution, in this case.
\end{rmk}

Collecting the results of the last two subsections, we can describe an equilibrium of the proposed network game.

\begin{theorem}\label{thm:game.main}
Assume that, for every $x\in\mX$, $Y^x_0\sim p_0(y;x)\,\mathrm{d}y$ with $p_0(\cdot;x)\in W^1_2(\RR_+)$ and $p_0(0;x)=0$.
Denote by $\mS$ the mapping constructed in Proposition \ref{prop:game.prop2}, with the associated $(\hat{r},\hC^b,\hC^l,\hnu)$, and by $\blambda$ the unique fixed-point of the associated mapping $\mD$. Then, $(r,\theta,\bC^b,\bC^l,\bnu)$, with
\begin{equation*}
\begin{split}
r_t(x):=\hat{r}(\blambda_t;x), \,\theta_t(x):= \hat{\theta}_t(\blambda;x),
\,\bC^b_t(x):=\hC^b(\blambda_t;x), \\
\bC^l_t(x):=\hC^l(\blambda_t;x),
\, \bnu_t(x,x'):=\hnu(\blambda_t;x,x'),
\end{split}
\end{equation*}
form an equilibrium.
\end{theorem}

\noindent\textbf{Proof.} The theorem follows by collecting the established results. First, we consider $\mS$ from Proposition \ref{prop:game.prop2} and use Lemma \ref{le:game.staticSuffCond}, equation (\ref{eq.game.b.Static}), as well as the Lipschitz and boundedness properties of the associated $R_t$, to conclude that $\mS$ satisfies Assumption \ref{ass:game.1}. Finally, we consider the fixed-point $\blambda$ from Proposition \ref{prop:game.prop1} and use Corollary \ref{cor:game.EquilSuffCond.Dyn.2} to conclude that $(r,\theta,\bC^b,\bC^l,\bnu)$ form an equilibrium. \ep

\medskip

It is important to note that the results of this section are not limited to proving the existence of an equilibrium: the existence proof provided herein is, in fact, constructive, and it also produces a computational algorithm. Indeed, the proof of Theorem \ref{thm:game.main} shows that the equilibrium quantities $(r,\theta,\bC^b,\bC^l,\bnu)$ are computed through the associated $\mS$ and $\blambda$. The mapping $\mS$ can be easily implemented (i.e. programmed) via the algorithm described in the proof of Proposition \ref{prop:game.prop2}. Finally, the fixed-point $\blambda$ can be obtained by the Picard iteration described in the proof of Proposition \ref{prop:game.prop1}, which, in turn, can be computed numerically using $\mS$ and an approximation scheme for the associated linear parabolic PDEs.

\section{Appendix A}

\noindent\textbf{Proof of Theorem \ref{thm existence}. Step 1.} Notice that the theorem amounts to the existence of a family $\{\Lambda^x\}_{x\in\mathcal{X}}$ of (``loss'') processes satisfying the fixed-point constraints
\begin{equation}\label{def:iteration}
\Lambda^x_t=C(x)\,\int_{\mathcal{X}} g\Big(\pp\Big(\min_{s\in[0,t]} \big(Z^{x'}_s+\Lambda^{x'}_s\big)>0\Big)\!\Big)\,\kappa(x,\mathrm{d}x'),\;\;t\in[0,T],\quad x\in\mathcal{X}
\end{equation}
and converging to zero as $t\downarrow0$.
We aim to employ Theorem \ref{thm:schauder} and, to this end, start with the construction of a domain $\mathscr{D}$ as therein. 

\medskip

Notice that Assumption \ref{asmp_ex}(a) and the continuity of $Z^x$, $x\in\mX$ yield the existence of a family $\{(\delta^x_n,\eta^x_n)_{n\in\nn}\}_{x\in\mathcal{X}}$ of positive sequences, satisfying $\delta^x_n,\eta^x_n\downarrow0$, as $n\to\infty$, for all $x\in\mathcal{X}$, and
\begin{equation}\label{deltan_choice}
\int_{\mathcal X} g\Big(\pp\Big(\min_{s\in[0,\delta^x_n]} Z^{x'}_s>\eta^{x'}_n\Big)\!\Big)\,\kappa(x,\mathrm{d}x')>-\frac{\eta^x_n}{C(x)},\quad x\in\mathcal{X},\quad n\in\nn.
\end{equation}
We define $\mathscr{D}\subset D([0,T],\rr)^{\mathcal{X}}$ as the set of families $\{\Lambda^x\}_{x\in\mathcal X}$ of c\`{a}dl\`{a}g non-increasing functions on $[0,T]$, with $\Lambda^x_0=0$, for all $x\in\mathcal{X}$, and
\begin{equation}\label{def:domain}
\Lambda^x_t\ge -\eta^x_n,\;\; t\in[\delta^x_{n+1},\delta^x_n),\quad n\in\nn,\quad x\in\mathcal{X}.
\end{equation}

\smallskip

Clearly, $\mathscr{D}\subset D([0,T],\rr)^{\mathcal{X}}$ is convex. In addition, since every convergent sequence in $D([0,T],\rr)$ converges pointwise on a dense subset of $[0,T]$ which includes $0$ (see e.g. \cite[Theorem 12.5.1(iv)]{Whitt}) and each function $t\mapsto\sum_{n=1}^\infty -\eta^x_n\,\mathbf{1}_{[\delta^x_{n+1},\delta^x_n)}(t)$ is right-continuous on $[0,\delta^x_1)$, the set $\mathscr{D}$ is also closed. Moreover, $S([0,T],\rr)^{\mathcal{X}}\cap\mathscr{D}$ is dense in $\mathscr{D}$, as can be inferred from the definition of the strong M1-metric in terms of graph parametrizations (see e.g. \cite[Section 12.3]{Whitt}).

\medskip

\noindent\textbf{Step 2.} Next, we factor the map associated with the right-hand side of \eqref{def:iteration} into a composition $F_2\circ F_1$ of two auxiliary maps $F_1$ and $F_2$. The components of the map
\begin{equation}
F_1:\;\mathscr{D}\longrightarrow \big({\mathscr P}\big(D([0,T],\rr)^2\big)\big)^{\mathcal{X}}
\end{equation} 
send each family $\{\Lambda^x\}_{x\in\mathcal{X}}\in\mathscr{D}$ to the joint law of the pair $(Z^x,\tilde{\Lambda}^x)$, where
\begin{equation}\label{def:tilde_lambda}
\tilde{\Lambda}^x_t:=\begin{cases}
\Lambda^x_t\qquad\,\text{if}\;\;\;t\in[0,\tau^x\wedge T), \\
-Z^x_{\tau^x}\;\;\;\text{if}\;\;\;t\in(\tau^x\wedge T,T].
\end{cases}
\end{equation}
Here, ${\mathscr P}(D([0,T],\rr)^2)$ is the space of Borel probability measures on $D([0,T],\rr)^2$ with the topology of weak convergence, and $({\mathscr P}(D([0,T],\rr)^2))^{\mathcal{X}}$ is the corresponding product space with the product topology. The two key properties of $F_1$ are 
\begin{eqnarray}
&& \pp\Big(\min_{s\in[0,t]} \big(Z^x_s+\tilde{\Lambda}^x_s\big)>0\Big)
=\pp\Big(\min_{s\in[0,t]} \big(Z^x_s+\Lambda^x_s\big)>0\Big),\;\;t\in[0,T],\quad x\in\mathcal{X}, \label{Lambda_tilde_Lambda} \\
&& Z^x_t+\tilde{\Lambda}^x_t=0,\;\;t\in(\tau^x\wedge T,T],\quad x\in\mathcal{X}.
\end{eqnarray}
In particular, with 
\begin{equation}
\begin{split}
& F_2:\;\mathrm{ran\;}F_1\longrightarrow D([0,T],\rr)^{\mathcal{X}}, \\
& \{\mathscr{L}(Z^x,\tilde{\Lambda}^x)\}_{x\in\mathcal{X}}\mapsto \bigg\{C(x)\int_{\mathcal{X}} g\Big(\pp\Big(\min_{s\in[0,t]} \big(Z^{x'}_s+\tilde{\Lambda}^{x'}_s\big)>0\Big)\!\Big)\,\kappa(x,\mathrm{d}x'),\;t\in[0,T]\bigg\}_{x\in\mathcal{X}},
\end{split}
\end{equation}
the composition $F_2\circ F_1$ is given by the right-hand side of \eqref{def:iteration}. 

\medskip

\noindent\textbf{Step 3.} The goal of this step is to show that the range of $F_1$ in $({\mathscr P}(D([0,T],\rr)^2))^{\mathcal{X}}$ is compact. Since each component $F^x_1$ of $F_1$ is defined in terms of $\Lambda^x$ and $Z^x$ only, the range of $F_1$ is of product form and, by Tychonoff's theorem, it is enough to verify the compactness of the range of each $F^x_1$ in ${\mathscr P}(D([0,T],\rr)^2)$. By Prokhorov's theorem, the latter is precompact if, for any $\iota\in(0,1)$, one can find a compact subset of $D([0,T],\rr)^2$ which is assigned a probability of at least $1-\iota$ by all the probability measures in the range of $F^x_1$. In view of the compactness criterion for subsets of $D([0,T],\rr)$ (cf. \cite[Theorem 12.12.2]{Whitt}), such is possible thanks to the first marginal of $F_1^x(\{\Lambda^{x'}\}_{x'\in\mathcal X})$ being the law of $Z^x$ for all $\{\Lambda^{x'}\}_{x'\in\mathcal{X}}\in\mathscr{D}$, the inequalities
\begin{equation}
0\ge\tilde{\Lambda}^x_t\ge\min_{s\in[0,T]} (-Z^x_s),\quad t\in[0,T]
\end{equation}
(see the definition of $\tilde{\Lambda}^x$ in \eqref{def:tilde_lambda}), the estimates
\begin{equation}
\tilde{\Lambda}^x_t=\Lambda^x_t\ge -\eta^x_n,\;\;t\in[0,\delta^x_n)\quad\text{on}\quad\{\tau^x\ge\delta^x_n\}
\end{equation}
for all $n\in\nn$ (cf. \eqref{def:domain}), and the monotonicity of $\tilde{\Lambda}^x$.  

\medskip

It remains to prove that the range of $F^x_1$ is closed in ${\mathscr P}(D([0,T],\rr)^2)$ and, thus, compact. Consider a sequence $\mathcal{L}(Z^x,\tilde{\Lambda}^{x,k})$, $k\in\nn$ in the range of $F_1^x$ converging weakly to some $\mathcal{L}(Z^x,\tilde{\Lambda}^{x,\infty})\in{\mathscr P}(D([0,T],\rr)^2)$. Due to the Skorokhod representation theorem in the form of \cite[Theorem 3.5.1]{Dud}, there exist $(\underline{Z}^{x,k},\underline{\tilde{\Lambda}}^{x,k})\stackrel{d}{=}(Z^x,\tilde{\Lambda}^{x,k})$, $k\in\nn$ and $(\underline{Z}^{x,\infty},\underline{\tilde{\Lambda}}^{x,\infty})\stackrel{d}{=}(Z^x,\tilde{\Lambda}^{x,\infty})$ on some probability space $(\Omega,\mathcal{F},\pp)$ such that $(\underline{Z}^{x,k},\underline{\tilde{\Lambda}}^{x,k})\to (\underline{Z}^{x,\infty},\underline{\tilde{\Lambda}}^{x,\infty})$ as $k\to\infty$
on $\Omega\backslash\Omega_0$ for some $\pp$-null set $\Omega_0$. For any  $\omega\in\Omega\backslash\Omega_0$, the sequence 
\begin{equation}
\underline{\tau}^{x,k}(\omega):=\inf\{t\in[0,T]:\,\underline{Z}^{x,k}(\omega)+\underline{\tilde{\Lambda}}^{x,k}(\omega)=0\}\wedge T,\quad k\in\nn
\end{equation} 
has a subsequence converging to a limit $\underline{\tau}^{x,\infty}(\omega)\in[0,T]$. In the rest of this step, we fix an $\omega\in\Omega\backslash\Omega_0$ and take all limits $k\to\infty$ along the  subsequence associated with that $\omega$. 

\medskip

One has $\underline{\tilde{\Lambda}}^{x,k}_t(\omega)\to \underline{\tilde{\Lambda}}^{x,\infty}_t(\omega)$ on a dense subset of $[0,T]$ that includes $0$, which together with the definition of $\tilde{\Lambda}^x$ in \eqref{def:tilde_lambda} implies 
\begin{equation}
\underline{\tilde{\Lambda}}^{x,\infty}_t(\omega)=\begin{cases}
\Lambda^{x,\infty}_t,\;\;t\in[0,\underline{\tau}^{x,\infty}(\omega)), \\
-\lim_{k\to\infty}\,\underline{Z}^{x,k}_{\underline{\tau}^{x,k}(\omega)}(\omega),\;\;t\in(\underline{\tau}^{x,\infty}(\omega),T],
\end{cases}
\end{equation} 
where $\Lambda^{x,\infty}$ is the M1-limit of $\Lambda^{x,k}$ on a suitable sequence of closed intervals increasing to $[0,\underline{\tau}^{x,\infty}(\omega))$, and $\Lambda^{x,k}$ belongs to the preimage of $\mathcal{L}(Z^x,\tilde{\Lambda}^{x,k})$ under $F^x_1$. The closedness assertion now follows from
\begin{equation}
\begin{split}
\underline{Z}^{x,\infty}_{\underline{\tau}^{x,\infty}(\omega)}(\omega)+\underline{\tilde{\Lambda}}^{x,\infty}_{\underline{\tau}^{x,\infty}(\omega)}(\omega)
& =\lim_{t\downarrow\underline{\tau}^{x,\infty}(\omega)}\,\lim_{k\to\infty}\, \big(\underline{Z}^{x,k}_t(\omega)+\underline{\tilde{\Lambda}}^{x,k}_t(\omega)\big) \\
& =\lim_{t\downarrow\underline{\tau}^{x,\infty}(\omega)}\,\lim_{k\to\infty}\, \big(\underline{Z}^{x,k}_t(\omega)-\underline{Z}^{x,k}_{\underline{\tau}^{x,k}(\omega)}(\omega)\big)=0,
\end{split}
\end{equation}
with the last equality being a result of the uniform convergence $\underline{Z}^{x,k}(\omega)\stackrel{k\to\infty}{\longrightarrow}\underline{Z}^{x,\infty}(\omega)$ (see e.g. \cite[Theorem 4.2]{Delarue2}).

\medskip

\noindent\textbf{Step 4.} Next, we argue that the map $F_2$ is continuous. Combined with the step 3 this yields the compactness of the range of $F_2$. To deduce the continuity of $F_2$ we pick a sequence
$\{{\mathcal L}(Z^x,\tilde{\Lambda}^{x,k})\}_{x\in\mathcal X}=F_1(\{\Lambda^{x,k}\}_{x\in\mathcal{X}})$, $k\in\nn$ converging to some $\{{\mathcal L}(Z^x,\tilde{\Lambda}^{x,\infty})\}_{x\in\mathcal X}=F_1(\{\Lambda^{x,\infty}\}_{x\in\mathcal{X}})$. For any $x\in\mathcal{X}$ and any continuity point $t\in[0,T]$ of $\Lambda^{x,\infty}$, the process $\tilde{\Lambda}^{x,\infty}$ is continuous at $t$, so that one can infer $\lim_{k\to\infty} \mathcal{L}(Z^x,\tilde{\Lambda}^{x,k})=\mathcal{L}(Z^x,\tilde{\Lambda}^{x,\infty})$ in ${\mathscr P}(D([0,t],\rr)^2)$ from the Skorokhod representation theorem and the characterizations of M1-convergence for continuous and monotone functions (see e.g. \cite[Theorem 4.2]{Delarue2}). Hence, $\lim_{k\to\infty} \mathcal{L}(Z^x+\tilde{\Lambda}^{x,k})=\mathcal{L}(Z^x+\tilde{\Lambda}^{x,\infty})$ in ${\mathscr P}(D([0,t],\rr))$ (see e.g. \cite[Theorem 12.7.3]{Whitt}). A subsequent application of the Skorokhod representation theorem and the definition of the strong M1-metric in terms of graph parametrizations gives the convergences in distribution
\begin{equation}\label{eq:min_conv}
\lim_{k\to\infty}\,\min_{s\in[0,t]} \big(Z^x_s+\tilde{\Lambda}^{x,k}_s\big) =\min_{s\in[0,t]} \big(Z^x_s+\tilde{\Lambda}^{x,\infty}_s\big)
\end{equation} 
for all continuity points $t\in[0,T]$ of $\Lambda^{x,\infty}$ and all $x\in\mathcal{X}$. The latter argument shows \eqref{eq:min_conv} also for $t=T$ and all $x\in\mathcal{X}$. 

\medskip

Further, we claim that, for any $x\in\mathcal{X}$, 
\begin{equation}\label{eq:log_conv}
\lim_{k\to\infty}\,g\Big(\pp\Big(\min_{s\in[0,t]} \big(Z^x_s+\tilde{\Lambda}^{x,k}_s\big)>0\Big)\!\Big)=g\Big(\pp\Big(\min_{s\in[0,t]} \big(Z^x_s+\tilde{\Lambda}^{x,\infty}_s\big)>0\Big)\!\Big)
\end{equation} 
at all continuity points $t\in[0,T]$ of $\Lambda^{x,\infty}$ and at $t=T$. By \eqref{eq:min_conv} the convergence of the probabilities in \eqref{eq:log_conv} can be reduced to
\begin{equation}\label{no_atom_identity}
\pp\Big(\min_{s\in[0,t]} \big(Z^x_s+\tilde{\Lambda}^{x,\infty}_s\big)=0\Big)=0.
\end{equation}
The identity \eqref{no_atom_identity} holds, since the intersection of the events $\min_{s\in[0,t]} (Z^x_s+\tilde{\Lambda}^{x,\infty}_s)=0$ and
\begin{equation}
\inf\{s\in[0,t]:\,Z^x_s+\tilde{\Lambda}^{x,\infty}_s\le 0\}=t 
\end{equation} 
is contained in the zero probability event $Z^x_t=-\Lambda^{x,\infty}_t$ (see Assumption \ref{asmp_ex}(b)), and the intersection of the events $\min_{s\in[0,t]} (Z^x_s+\tilde{\Lambda}^{x,\infty}_s)=0$ and
\begin{equation}
\inf\{s\in[0,t]:\,Z^x_s+\tilde{\Lambda}^{x,\infty}_s\le 0\}<t 
\end{equation} 
has a probability of zero by Assumption \ref{asmp_ex}(c). At this point, \eqref{eq:log_conv} is readily obtained from the convergence of the probabilities involved and the positivity of the probability on the right-hand side of \eqref{eq:log_conv} (recall \eqref{Lambda_tilde_Lambda} and Assumption \ref{asmp_ex}(d)).


\medskip

\noindent\textbf{Step 5.} Lastly, for any sequence $\{\Lambda^{x,k}\}_{x\in\mathcal{X}}$, $k\in\nn$ in ${\mathscr D}$ converging to some $\{\Lambda^{x,\infty}\}_{x\in\mathcal{X}}\in{\mathscr D}$, we proceed as in the step 4
to conclude that $(Z^x+\Lambda^{x,k})\stackrel{k\to\infty}{\longrightarrow}Z^x+\Lambda^{x,\infty}$ in $D([0,t],\rr)$ almost surely for all $x\in\mathcal{X}$ and all continuity points $t\in[0,T]$ of $\Lambda^{x,\infty}$ and $t=T$. As in the step 4, this leads to 
\begin{eqnarray}
&& \lim_{k\to\infty}\,\min_{s\in[0,t]} \big(Z^x_s+\Lambda^{x,k}_s\big)=\min_{s\in[0,t]}\,\big(Z^x_s+\Lambda^{x,\infty}_s\big)\;\;\text{almost surely}, \\
&& \lim_{k\to\infty}\,g\Big(\pp\Big(\min_{s\in[0,t]} \big(Z^x_s+\Lambda^{x,k}_s\big)>0\Big)\!\Big)=g\Big(\pp\Big(\min_{s\in[0,t]} \big(Z^x_s+\Lambda^{x,\infty}_s\big)>0\Big)\!\Big),
\end{eqnarray}
for all $t$ and $x$ as described. We can therefore infer from the M1-convergence criterion for monotone functions in \cite[Theorem 4.2]{Delarue2} that the map
\begin{equation}
\begin{split}
& F_2\circ F_1:\;\mathscr{D}\longrightarrow D([0,T],\rr)^{\mathcal{X}}, \\
& \{\Lambda^x\}_{x\in\mathcal{X}}\mapsto \bigg\{C(x)\int_{\mathcal{X}} g\Big(\pp\Big(\min_{s\in[0,t]} \big(Z^{x'}_s+\Lambda^{x'}_s\big)>0\Big)\!\Big)\,\kappa(x,\mathrm{d}x'),\;t\in[0,T]\bigg\}_{x\in\mathcal{X}},
\end{split}
\end{equation}
is continuous. 

\medskip

In addition, for any $\{\Lambda^x\}_{x\in\mathcal{X}}\in\mathscr{D}$, it holds
\begin{equation}
\begin{split}
& \;C(x)\,\int_{\mathcal{X}} g\Big(\pp\Big(\min_{s\in[0,t]} \big(Z^{x'}_s+\Lambda^{x'}_s\big)>0\Big)\!\Big)\,\kappa(x,\mathrm{d}x') \\
& \ge C(x)\,\int_{\mathcal{X}} g\Big(\pp\Big(\min_{s\in[0,t]} Z^{x'}_s>\eta^{x'}_n\Big)\!\Big)\,\kappa(x,\mathrm{d}x')
>-\eta^x_n,\quad t\in[\delta^x_{n+1},\delta_n^x),\quad n\in\nn
\end{split}
\end{equation}
thanks to \eqref{deltan_choice}, so that the compact $\mathrm{ran}\,(F_2\circ F_1)=\mathrm{ran}\,F_2$ (see the step 4) is a subset of $\mathscr{D}$. In summary, $F_2\circ F_1$ maps the convex closed set $\mathscr{D}$, in which $S([0,T],\rr)^{\mathcal{X}}\cap\mathscr{D}$ forms a dense subset, continuously onto a compact subset of $\mathscr{D}$. Thus, $F_2\circ F_1$ has a fixed-point by Theorem \ref{thm:schauder}. \ep

\section{Appendix B}

\begin{lemma}\label{thm:jump_crit}
Fix an arbitrary $t\in[0,T)$ and suppose that Assumptions \ref{asmp_ex}(b) and \ref{asmp_dom}(a) hold. Assume, in addition, that a solution $\{Y^x\}_{x\in\mathcal{X}}$ of \eqref{what is sol} admits $\{k^x\ge0\}_{x\in\mathcal X}$ such that, for each $x\in\mathcal X$,
\begin{equation}\label{loc_slope}
g(\pp(\tau^x>t)\!)-g(z)\ge k^x\big(\pp(\tau^x>t)-z\big)
\end{equation}
for all $z$ in a left neighborhood of $\pp(\tau^x>t)$, and $\{c^x\}_{x\in\mX}$, $\{z^x\}_{x\in\mX}$ satisfying \eqref{eq.cz.def} and such that
\begin{equation}\label{eq:jump series}
\begin{split}
\sum_{n=1}^\infty \,\int_{\mathcal{X}^n} \prod_{m=0}^{n-1} \bigg(C(x_m)\frac{c^{x_{m+1}}\,k^{x_{m+1}}}{g'(\pp(\tau^{x_{m+1}}\!>\!t)\!)}\bigg) \prod_{m=1}^{n-1} 
\pp(\xi^{x_m}_{s_{m+1},s_m}\!\le\! z^{x_m}\!-\!\eps^{x_m})\,\E[\xi^{x_n}_{s_{n+1},s_n}\!\wedge\! z^{x_n}]\qquad\;\; \\ \kappa(x_0,\mathrm{d}x_1)\,\kappa(x_1,\mathrm{d}x_2)\,\ldots\,\kappa(x_{n-1},\mathrm{d}x_n)=\infty,
\end{split}
\end{equation}
for some $x_0\in\mathcal{X}$, $(s_n)_{n\in\nn}\downarrow t$, and $\eps^x\in(0,z^x)$, $x\in\mathcal{X}$. Then, $Y^x_{t-}>Y^x_t$ for at least one $x\in\mathcal{X}$. 
\end{lemma}

\noindent\textbf{Proof.} Assume the opposite, i.e. that there exist $t$, $\{Y^x\}_{x\in\mathcal{X}}$, $\{k^x\}_{x\in\mX}$, $\{c^x\}_{x\in\mX}$, $\{z^x\}_{x\in\mX}$, $x_0$, $(s_n)_{n\in\nn}$ and $\{\eps^x\}_{x\in\mX}$ satisfying the assumptions of the lemma and such that $Y^x_{t-}=Y^x_t$ for all $x\in\mathcal{X}$. Then, Assumption \ref{asmp_ex}(b) gives 
\begin{equation}\label{>t bound}
\begin{split}
\pp\big(\tau^x>t,\,Y^x_t\in(0,z)\!\big)\,g'(\pp(\tau^x>t)\!)
=\pp\big(\tau^x\ge t,\,Y^x_{t-}\in(0,z)\!\big)\,g'(\pp(\tau^x\ge t)\!)\ge c^x z,  \\
z\in[0,z^x],\quad x\in\mathcal{X}.
\end{split}
\end{equation}

At this point, we let
\begin{eqnarray}
&& \Lambda^x_s:=C(x)\int_{\mathcal X} g(\pp(\tau^{x'}>s)\!)\,\kappa(x,\mathrm{d}x'),\;\; s\in[0,T],\quad x\in\mathcal{X}, \\
&& \zeta^x(\cdot):=\frac{\pp(\tau^x>t,\,Y^x_t\in\cdot\,)}{\pp(\tau^x>t)},\quad x\in\mathcal{X}.
\end{eqnarray}
We assume without loss of generality that 
\begin{equation}\label{eps_bound}
\Lambda^x_t-\Lambda^x_{s_1}\le \eps^x,\quad x\in\mathcal{X},
\end{equation}
since \eqref{eq:jump series} remains valid for any $(s_n)_{n\ge N}$ in place of $(s_n)_{n\in\nn}$.
Next, we employ \eqref{loc_slope}, the monotonicity of $\{\Lambda^x\}_{x\in\mathcal{X}}$, the inequality \eqref{eq:dom}, the stochastic dominance in \eqref{>t bound}, the formula $\int_0^\infty \pp(R\ge z)\,\mathrm{d}z=\E[R]$ for non-negative random variables $R$, and the elementary bound $(z_1+z_2)\wedge z_3\ge z_1\wedge z_3+z_2\,\mathbf{1}_{\{z_1+z_2\le z_3\}}$, $z_1,z_2,z_3\in\rr$ in conjunction with \eqref{eps_bound}, consecutively, to obtain
\begin{equation}\label{jump_basic_est}
\begin{split}
&\,\Lambda^{x_0}_{s_1}-\Lambda^{x_0}_t \le  C(x_0)\int_{\mathcal X} 
k^{x_1}\,\pp(\tau^{x_1}>t) \\
& \qquad\qquad\qquad\qquad
\cdot\!\!\Big(\pp\Big(\min_{s\in[t,s_1]}\!\big(Y^{x_1}_t\!+\!Z^{x_1}_s\!-\!Z^{x_1}_t\!+\!\Lambda^{x_1}_s\!-\!
\Lambda^{x_1}_t\big)\!>\!0\,\Big|\,\tau^{x_1}\!>\!t\Big)\!-\!1\Big)\,\kappa(x_0,\mathrm{d}x_1) \\
&\le -C(x_0)\int_{\mathcal X}  \int_0^{z^{x_1}} \pp(\xi^{x_1}_{s_2,s_1}\ge z+\Lambda^{x_1}_{s_2}-\Lambda^{x_1}_t)\,\zeta^{x_1}(\mathrm{d}z)\,\pp(\tau^{x_1}>t)\,k^{x_1}\,\kappa(x_0,\mathrm{d}x_1) \\
& \le -C(x_0)\int_{\mathcal X} \frac{c^{x_1}\,k^{x_1}}{g'(\pp(\tau^{x_1}>t)\!)}\,\int_0^{z^{x_1}} \pp(\xi^{x_1}_{s_2,s_1}\ge z+\Lambda^{x_1}_{s_2}-\Lambda^{x_1}_t)\,\mathrm{d}z\,\kappa(x_0,\mathrm{d}x_1) \\
&= -C(x_0)\int_{\mathcal X}  \frac{c^{x_1}\,k^{x_1}}{g'(\pp(\tau^{x_1}>t)\!)}\,\E\big[(\xi^{x_1}_{s_2,s_1}+\Lambda^{x_1}_t-\Lambda^{x_1}_{s_2})\wedge z^{x_1}\big]\,\kappa(x_0,\mathrm{d}x_1) \\
&\le -\,C(x_0)\int_{\mathcal X} \frac{c^{x_1}\,k^{x_1}}{g'(\pp(\tau^{x_1}>t)\!)}\,\E[\xi^{x_1}_{s_2,s_1}\wedge z^{x_1}]\,\kappa(x_0,\mathrm{d}x_1) \\
&\quad \,+C(x_0)\int_{\mathcal X} \frac{c^{x_1}\,k^{x_1}}{g'(\pp(\tau^{x_1}>t)\!)}\,(\Lambda^{x_1}_{s_2}-\Lambda^{x_1}_t)\,\pp(\xi^{x_1}_{s_2,s_1}\le z^{x_1}-\eps^{x_1})\,\kappa(x_0,\mathrm{d}x_1).
\end{split}
\end{equation}

To conclude we repeat the estimates of \eqref{jump_basic_est} mutatis mutandis for the term $\Lambda^{x_1}_{s_2}-\Lambda^{x_1}_t$ inside the last expression in \eqref{jump_basic_est} to find
\begin{equation}
\begin{split}
&\Lambda^{x_0}_{s_1}-\Lambda^{x_0}_t\le -C(x_0) \int_{\mathcal X} \frac{c^{x_1}\,k^{x_1}}{g'(\pp(\tau^{x_1}>t)\!)}\,\E[\xi^{x_1}_{s_2,s_1}\wedge z^{x_1}]\,\kappa(x_0,\mathrm{d}x_1) \\
&\qquad\qquad\quad\;\, -C(x_0) \int_{\mathcal{X}} \frac{c^{x_1}\,k^{x_1}}{g'(\pp(\tau^{x_1}>t)\!)}\,\pp(\xi^{x_1}_{s_2,s_1}\le z^{x_1}-\eps^{x_1})\,\\
&\qquad\qquad\qquad\qquad\quad\;\;\;
\cdot C(x_1)\int_{\mathcal{X}} \frac{c^{x_2}\,k^{x_2}}{g'(\pp(\tau^{x_2}>t)\!)}\,\E[\xi^{x_2}_{s_3,s_2}\wedge z^{x_2}]\,\kappa(x_1,\mathrm{d}x_2)\,\kappa(x_0,\mathrm{d}x_1) \\
&\qquad\qquad\quad\;\, +C(x_0)\int_{\mathcal{X}} \frac{c^{x_1}\,k^{x_1}}{g'(\pp(\tau^{x_1}>t)\!)}\,\pp(\xi^{x_1}_{s_2,s_1}\le z^{x_1}-\eps^{x_1})\\
&\qquad\qquad\qquad\qquad\qquad\, \cdot C(x_1)\int_{\mathcal{X}} \frac{c^{x_2}\,k^{x_2}}{g'(\pp(\tau^{x_2}>t)\!)}\,(\Lambda^{x_2}_{s_3}-\Lambda^{x_2}_t)\,\pp(\xi^{x_2}_{s_3,s_2}\le z^{x_2}-\eps^{x_2}) \\ 
&\qquad\qquad\qquad\qquad\qquad\qquad\qquad\qquad\qquad\qquad \qquad\qquad\qquad\; 
\kappa(x_1,\mathrm{d}x_2) \, \kappa(x_0,\mathrm{d}x_1).
\end{split}
\end{equation} 
Iterating $N\in\nn$ times and dropping the term involving $\Lambda^{x_N}_{s_{N+1}}-\Lambda^{x_N}_t\le 0$ we arrive at 
\begin{equation}
\begin{split}
& \Lambda^{x_0}_{s_1}-\Lambda^{x_0}_t \\ 
&\le -\sum_{n=1}^N \int_{\mathcal{X}^n} \prod_{m=0}^{n-1} \bigg(\!C(x_m)\frac{c^{x_{m+1}}\,k^{x_{m+1}}}{g'(\pp(\tau^{x_{m+1}}\!>\!t)\!)}\!\bigg) \prod_{m=1}^{n-1}  \pp(\xi^{x_m}_{s_{m+1},s_m}\le z^{x_m}\!-\!\eps^{x_m})\,\E[\xi^{x_n}_{s_{n+1},s_n}\wedge z^{x_n}] \\ 
& \qquad\qquad\qquad\qquad\qquad\qquad\qquad\qquad\qquad\qquad\quad
\kappa(x_0,\mathrm{d}x_1)\,\kappa(x_1,\mathrm{d}x_2)\,\ldots\,\kappa(x_{n-1},\mathrm{d}x_n).
\end{split}
\end{equation}  
We now take $N\to\infty$ and end up with $\Lambda^{x_0}_{s_1}-\Lambda^{x_0}_t\le  -\infty$ (recall \eqref{eq:jump series}), which is the desired contradiction. \ep

\medskip

Using the above lemma and making a somewhat stronger assumption on $Z$, we obtain a sufficient condition for fragility that is, in some sense, sharper than Theorem \ref{cor_PF} (i.e. it covers the case $\varrho=0$).

\begin{theorem}\label{cor_semimart}
Fix an arbitrary $t\in[0,T)$ and suppose that, for each $x\in\mathcal{X}$, 
\begin{equation}
g(\pp(\tau^x>t)\!)-g(z)\ge g'(\pp(\tau^x>t)\!)\big(\pp(\tau^x>t)-z\big)
\end{equation}
for all $z$ in a left neighborhood of $\pp(\tau^x>t)$ and that $\{Z^x\}_{x\in\mX}$ satisfy (\ref{semimart_set}), with the associated $\overline{\alpha}\in\rr$ and $\underline{\sigma},\overline{\sigma}\in(0,\infty)$.
Assume that a solution $\{Y^x\}_{x\in\mathcal{X}}$ of \eqref{what is sol} admits $\{c^x\}_{x\in\mX}$, $\{z^x\}_{x\in\mX}$ satisfying \eqref{eq.cz.def} and such that
\begin{equation}\label{semimart_cond}
\begin{split}
& \sum_{n=1}^\infty \big(\sqrt{s_n-s_{n+1}}+\tilde{\alpha}(s_{n+1}-t)\big) \\
& \quad\;\;\cdot\int_{\mathcal{X}^n} \prod_{m=0}^{n-1} (C(x_m)c^{x_{m+1}}) 
\prod_{m=1}^{n-1}
\bigg(1-4\overline{\Phi}\bigg(\frac{z^{x_m}-\eps^{x_m}+\overline{\alpha}(s_m-t)}
{\sqrt{\underline{\sigma}^2(s_m\!-\!s_{m+1})+\overline{\sigma}^2(s_{m+1}\!-\!t)}}\bigg)
\!\bigg) \\
&\qquad\qquad\qquad\qquad\qquad\qquad\qquad\qquad\;
\kappa(x_0,\mathrm{d}x_1)\,\kappa(x_1,\mathrm{d}x_2)\,\ldots\,\kappa(x_{n-1},\mathrm{d}x_n)=\infty
\end{split}
\end{equation}
for some $(s_n)_{n\in\nn}\downarrow t$, $\tilde{\alpha}<-\sqrt{\pi/2}\,\overline{\alpha}/\underline{\sigma}$,
$x_0\in\mathcal{X}$, and $\eps^x\in(0,z^x)$, $x\in\mX$. Then $Y^x_{t-}>Y^x_t$, for at least one $x\in\mathcal{X}$. Hereby, $\overline{\Phi}$ is the standard normal tail cumulative distribution function. 
\end{theorem}

\noindent\textbf{Proof.} We note that Assumption \ref{asmp_ex}(b) holds in the setting of \eqref{semimart_set} thanks to \cite[inequality (3)]{McNamara} in conjunction with Girsanov's theorem (see e.g. \cite[Chapter 3, Theorem 5.1 and Corollary 5.13]{KarShr}). 
Let us show that Assumption \ref{asmp_dom}(a) holds.
To this end, we write, for any $\theta\in(0,T-t)$, $t\le t_1<t_2\le t+\theta$, and $x\in\mathcal{X}$,
\begin{equation}
\min_{t_1\le s\le t_2} (Z^x_s-Z^x_t) = Z^x_{t_1}-Z^x_t+\min_{t_1\le s\le t_2} (Z^x_s-Z^x_{t_1}), 
\end{equation}
infer from \cite[inequality (1.5)]{Haj} that, conditionally on $Z^x_s$, $s\in[0,t]$,
\begin{equation}
Z^x_{t_1}-Z^x_t\le \overline{\alpha}(t_1-t)+\overline{\sigma}\frac{\beta^x_{t_1-t}+\tilde{\beta}^x_{t_1-t}}{2} 
\end{equation}
with (possibly dependent) standard Brownian motions $\beta^x$, $\tilde{\beta}^x$, and dominate the conditional law of $\min_{t_1\le s\le t_2} (Z^x_s-Z^x_{t_1})$ given $Z^x_{t_1}-Z^x_t$ and $Z^x_s$, $s\in[0,t]$ stochastically by $\overline{\alpha}(t_2-t_1)+\min_{0\le s\le\underline{\sigma}(t_2-t_1)} \hat{\beta}_s$ for a standard Brownian motion $\hat{\beta}$ by relying on the Dambis-Dubins-Schwarz theorem (see e.g. \cite[Chapter 3, Theorem 4.6]{KarShr}). Thus, Assumption \ref{asmp_dom}(a) holds with
\begin{equation}\label{xi_semimart}
\xi^x_{t_1,t_2}:=-\overline{\alpha}(t_2-t)-\overline{\sigma}\frac{\beta^x_{t_1-t}+\tilde{\beta}^x_{t_1-t}}{2}-\min_{0\le s\le\underline{\sigma}(t_2-t_1)} \hat{\beta}_s,
\end{equation} 
where $\hat{\beta}$ is a standard Brownian motion independent of all $\beta^x$, $\tilde{\beta}^x$.

\medskip

In view of the estimates 
\begin{equation}\label{semi_prob_est}
\begin{split}
\pp(\xi^x_{t_1,t_2}\le z^x\!-\!\eps^x) 
 =1-\pp\bigg(\!-\overline{\sigma}\frac{\beta^x_{t_1-t}\!+\!\tilde{\beta}^x_{t_1-t}}{2}+|\hat{\beta}_{\underline{\sigma}(t_2-t_1)}|
>z^x-\eps^x+\overline{\alpha}(t_2\!-\!t)\!\bigg) \\
 \ge 1-\frac{4}{\overline{\sigma}\sqrt{2\pi(t_1-t)}}\int_{-\infty}^\infty \overline{\Phi}\bigg(\frac{z^x-\eps^x+\overline{\alpha}(t_2-t)-z}{\underline{\sigma}\sqrt{t_2-t_1}}\bigg)\,e^{-\frac{z^2}{2\overline{\sigma}^2(t_1-t)}}\,\mathrm{d}z \\
 = 1-4\overline{\Phi}\bigg(\frac{z^x-\eps^x+\overline{\alpha}(t_2-t)}
{\sqrt{\underline{\sigma}^2(t_2-t_1)+\overline{\sigma}^2(t_1-t)}}\bigg)
\end{split}
\end{equation}
and
\begin{equation}\label{exp_est}
\begin{split}
& \E[\xi^x_{t_1,t_2}\wedge z^x]
 \ge \frac{1}{2}\,\E\Big[\big(\!-2\overline{\alpha}(t_1\!-\!t)\!-\!\overline{\sigma}(\beta^x_{t_1-t}\!+\!\tilde{\beta}^x_{t_1-t})\big)\!\wedge\! z^x\Big] \\
&\qquad\qquad\qquad\; +\frac{1}{2}\,\E\Big[\big(\!-2\overline{\alpha}(t_2\!-\!t_1)\!-\!2\min_{0\le s\le\underline{\sigma}(t_2-t_1)} \hat{\beta}_s\big)\!\wedge\! z^x\Big] \\
& = -\overline{\alpha}(t_2-t)+\underline{\sigma}\sqrt{\frac{2(t_2-t_1)}{\pi}} \\
&\quad\, +\frac{1}{2}\,\E\Big[\big(z^x+2\overline{\alpha}(t_1-t)+\overline{\sigma}(\beta^x_{t_1-t}+\tilde{\beta}^x_{t_1-t})\big)\mathbf{1}_{\{-2\overline{\alpha}(t_1-t)-\overline{\sigma}(\beta^x_{t_1-t}+\tilde{\beta}^x_{t_1-t})>z^x\}}\Big] \\
&\quad\,+\frac{1}{2}\,\E\Big[\big(z^x+2\overline{\alpha}(t_2-t_1)+2\min_{0\le s\le\underline{\sigma}(t_2-t_1)} \hat{\beta}_s\big)\mathbf{1}_{\{-2\overline{\alpha}(t_2-t_1)-2\min_{0\le s\le\underline{\sigma}(t_2-t_1)} \hat{\beta}_s>z^x\}}\Big] \\
& =  -\overline{\alpha}(t_1-t)+\underline{\sigma}\sqrt{\frac{2(t_2-t_1)}{\pi}}
+o(t_1-t)+o(\sqrt{t_2-t_1})\quad\text{as}\quad t_1,t_2\downarrow 0,
\end{split}
\end{equation}
the divergence of the series in \eqref{semimart_cond} implies the same for the series in \eqref{eq:jump series}, with the selection of $k^x=g'(\pp(\tau^x>t)\!)$ and of $\xi^x_{t_1,t_2}$ in \eqref{xi_semimart}. The result now follows from Lemma \ref{thm:jump_crit}.
In addition, (\ref{exp_est}) shows that Assumption \ref{asmp_dom}(b) holds with $s_n := t+\eta/n$, for sufficiently small $\eta>0$. \ep

\medskip

\noindent\textbf{Proof of Theorem \ref{cor_PF}. (a).} The continuity of $Z^x$, $x\in\mX$ and the defining equation \eqref{eq:dom} for $\xi^x_{t_1,t_2}$ imply that, for any $\{\varepsilon^x>0\}_{x\in\mX}$ and $(s_n)_{n\in\nn}\downarrow t$, we have
$$
\lim_{n\rightarrow\infty} \min_{x\in\mathcal{X}}\, \pp(\xi^x_{s_{n+1},s_n}\le \varepsilon^x) = 1,
$$
and, in turn,
\begin{equation}
\lim_{n\to\infty}\,\frac{1}{n}\log\prod_{m=1}^{n-1}  \big(\min_{x\in\mathcal{X}}\, \pp(\xi^x_{s_{m+1},s_m}\le \varepsilon^x)\big)=0. \label{regular1}
\end{equation}
In view of \eqref{what is rho}, \eqref{regular1} and \eqref{regular2}, for any $\{c^x, z^x\}_{x\in\mX}$ as in Theorem \ref{cor_PF}(a), any $\{\varepsilon^x\in(0,z^x)\}_{x\in\mX}$, and any $x_0$ that belongs to a closed irreducible component of $(C(x)\kappa(x,\{x'\})c^{x'})_{x,x'\in\mathcal{X}}$ with $\varrho>0$, we have:
\begin{equation}
\begin{split}
& \,\liminf_{n\to\infty}\,\frac{1}{n}\log\bigg(\int_{\mathcal{X}^n} \prod_{m=0}^{n-1} (C(x_m)c^{x_{m+1}}) \prod_{m=1}^{n-1} 
\pp(\xi^{x_m}_{s_{m+1},s_m}\le z^{x_m}-\eps^{x_m})\,\E[\xi^{x_n}_{s_{n+1},s_n}\wedge z^{x_n}] \\ 
& \qquad\qquad\qquad\qquad\qquad\qquad\qquad\qquad\qquad\qquad\qquad\quad\;\;\; \kappa(x_0,\mathrm{d}x_1)\,\ldots\,\kappa(x_{n-1},\mathrm{d}x_n)\bigg) \\
& \ge \liminf_{n\to\infty}\,\frac{1}{n}\log\bigg(\prod_{m=1}^{n-1} 
\big(\min_{x\in\mathcal{X}}\,\pp(\xi^x_{s_{m+1},s_m}\le z^x-\eps^x)\big) 
\,\min_{x\in\mathcal{X}}\,\E[\xi^x_{s_{n+1},s_n}\wedge z^x] \\
& \qquad\qquad\qquad\qquad\quad\;\;
\int_{\mathcal{X}^n} \prod_{m=0}^{n-1} (C(x_m)c^{x_{m+1}})\,   \kappa(x_0,\mathrm{d}x_1)\,\ldots\,\kappa(x_{n-1},\mathrm{d}x_n)\bigg)\ge\varrho>0. 
\end{split}
\end{equation}
Consequently, with $\{k^x\}_{x\in\mathcal X}$ obeying $\min_{x\in\mathcal X} (k^x/g'(\pp(\tau^x\!>\!t)\!)\!)>e^{-\rho}$, the summands in the series of \eqref{eq:jump series} grow exponentially and the result readily follows from Lemma \ref{thm:jump_crit}. 

\medskip

\noindent\textbf{(b).} The assumptions of Theorem \ref{cor_PF}(b) imply the existence of an eigenvector $v=(v^x)_{x\in\mX}$ of the matrix $A=(C(x)\kappa(x,\{x'\})c^{x'})_{x,x'\in\mathcal{X}}$ with strictly positive entries such that $Av < v$ (see e.g. \cite[Theorem 3.1.1(c)]{DZ}).
Then, $\eta^x_n:=v^x/n$ satisfy (\ref{near0_cond.t}) for all large enough $n$. Hence, Theorem \ref{thm existence} and Remark \ref{rem:ExistContSol} imply the existence of a solution to \eqref{what is sol}, with the prescribed distribution at $t-$ and continuous at $t$.
\ep

\bigskip\bigskip

\bibliographystyle{amsalpha}
\bibliography{MeanFieldNetworks}

\bigskip\bigskip

\end{document}